\documentclass{article}

\usepackage{a4wide}
\usepackage{authblk}
\usepackage{snapshot}
\date{}

\newcommand{\email}[1]{\texttt{#1},~}
\newcommand{\institution}[1]{#1\\}
\newcommand{\streetaddress}[1]{#1,~}
\newcommand{\city}[1]{#1,~}
\newcommand{\state}[1]{#1,~}
\newcommand{\postcode}[1]{#1,~}
\newcommand{\country}[1]{#1}

\usepackage{amsmath,amsfonts}
\usepackage{graphicx}
\usepackage{comment}
\usepackage{pdfsync}
\usepackage{todonotes}

\usepackage{hyperref}

\usepackage{tikz}
\usetikzlibrary{shapes,arrows,positioning,calc}

\tikzset{
  block/.style = {fill=white, minimum height=3em, minimum width=3em},
  core/.style={
           rectangle,
           rounded corners,
           draw=black, thin,
           minimum height=1.5em,
           inner sep=2pt,
           text centered,
         },
   optional/.style={
           rectangle,
           rounded corners,
           draw=black, thick,
           minimum height=1.5em,
           inner sep=2pt,
           text centered,
         },
}

\usepackage[draft,nomargin,inline,marginclue,multiuser]{fixme}

\FXRegisterAuthor{wb}{aom}{WB} 
\FXRegisterAuthor{lh}{alh}{LH}
\FXRegisterAuthor{mk}{amk}{MK}
\FXRegisterAuthor{am}{aam}{AM}

\newtheorem{remark}{Remark}
\newtheorem{primitive}{Primitive}

\newcommand{\dealii}{\textsc{deal.II}}
\newcommand{\aspect}{\textsc{Aspect}}

\begin{document}

\title{Propagating geometry information to finite element computations}

\author[1]{Luca Heltai}

\affil[1]{
 \email{luca.heltai@sissa.it}
 \institution{SISSA - International School for Advanced Studies}
 \streetaddress{mathLab, Mathematics Area, Via Bonomea 265}
 \city{34136 Trieste}
 \country{Italy}
}

\author[2]{Wolfgang Bangerth}
\affil[2]{
  \email{bangerth@colostate.edu}
  \institution{Colorado State University}
  \streetaddress{Department of Mathematics, Department of Geosciences, 1874 Campus Delivery}
  \city{Fort Collins}
  \state{CO}
  \postcode{80524}
  \country{USA}
}

\author[3]{Martin Kronbichler}
\affil[3]{
  \email{kronbichler@lnm.mw.tum.de}
  \institution{Technical University of Munich}
  \streetaddress{Institute for Computational Mechanics, Boltzmannstr.~15}
  \city{85748 Garching b. M\"unchen}
  \country{Germany}
}

\author[4]{Andrea Mola}
\affil[4]{
  \email{mola.andrea@sissa.it}
  \institution{SISSA - International School for Advanced Studies}
  \streetaddress{mathLab, Mathematics Area, Via Bonomea 265}
  \city{34136 Trieste}
  \country{Italy}
}






\maketitle

\section*{Abstract}
  The traditional workflow in continuum mechanics simulations is that a geometry description
  -- for example obtained using Constructive Solid Geometry or Computer Aided
  Design tools -- forms the input for a mesh generator. The mesh
  is then used as the sole input for the finite element, finite volume,
  and finite difference solver, which at this point no longer has
  access to the original, ``underlying'' geometry. However, many modern
  techniques -- for example, adaptive mesh refinement and the use
  of higher order geometry approximation methods -- really do need
  information about the underlying geometry to realize their full
  potential. We have undertaken an exhaustive study
  of where typical finite element codes use geometry information, with
  the goal of determining what information geometry tools would have
  to provide. Our study shows that nearly all geometry-related needs inside the
  simulators can be satisfied by just two ``primitives'': elementary
  queries posed by the simulation software to the geometry description.
  We then show that it is possible to provide these primitives in all
  of the commonly used ways in which geometries are described in
  common industrial workflows. We illustrate our solutions using
  examples from adaptive mesh refinement for complex geometries.
\section{Introduction}
\label{sec:introduction}

The traditional workflow of finite element, finite volume, and finite difference simulations
of physical processes consists of three phases: what is called ``preprocessing''; the actual 
numerical solution
of a partial differential equation; and what is called ``postprocessing''. In this workflow,
preprocessing generally means the generation of a geometric description of the domain
on which one wants to solve the problem -- either through the use of Computer Aided Design
(CAD), or by combining simpler geometries into one via constructive solid
geometry (CSG) -- and the
use of a mesh generator that uses the geometry to create the computational grid on which
the simulation is then run. On the other end of the pipeline, postprocessing consists of
the visualization of the computed solution and the extraction of quantities of interest.

Overall, this ``traditional'' workflow can be visualized through the following graph in which
information is only propagated from one box to the next:
\begin{center}
    \begin{tikzpicture}[auto, node distance=2cm,>=latex']
      \node [core, name=geometry] (geometry) {Geometry description};
      \node [core, right of=geometry, node distance=3.2cm] (mesh) {Mesh generation};
      \node [core, right of=mesh, node distance=3cm] (simulator) {Simulation};
      \node [core, right of=simulator, node distance=3cm] (visualization) {Visualization};
      \node [core, below=1em of visualization] (QoI) {Quantities of interest};

      \draw[->] (geometry) -- (mesh);
      \draw[->] (mesh) -- (simulator);
      \draw[->] (simulator) -- (visualization);
      \draw[->] (simulator) -- (QoI);
    \end{tikzpicture}
\end{center}
The fundamental issue with this workflow is that geometry information is only passed on
to the mesh generator, but is, in general, not available at the later stages of the pipeline.
This approach -- which to our knowledge is used in all commercial and open source 
simulation tools today -- is appropriate if the simulator is relatively simple; specifically, 
if (i) simulation and postprocessing tools only rely on a single, fixed mesh as their sole 
information on the geometry of the domain on which to solve the
problem under consideration, and (ii) if one uses lowest-order
finite-element, finite-volume, or finite-difference discretizations
for which it is sufficient to use the piecewise linear approximation
of the boundary that is generally obtained by replacing the ``underlying''
geometry of the problem by a fixed mesh characterized solely by its
vertices and assuming straight edges. In practice, the limitations of the workflow mentioned above
imply that to most finite element analysts, ``the
mesh is the domain'', even though to the designer the mesh is
generally an imperfect approximation of some underlying geometry
typically understood to be a CAD or CSG description.

Yet, simulation tools have become vastly more complex since the
formulation of the workflow above several decades ago, and as we will
show below, we can not make use of their full potential
\textit{unless geometry information is propagated to the simulation and analysis tools},
as well as to the postprocessing tools. For example, geometry information
is important in the following contexts:
\begin{itemize}
    \item Modern simulators no longer use only a single mesh, but create
      hierarchies of meshes. Two typical applications are the generation
      of a sequence of refined meshes to enable the use of geometric 
      multigrid solvers or preconditioners \cite{Briggs,Clevenger19}; and the use of
      adaptive mesh refinement to obtain a mesh
      that is better suited to the accurate solution of the underlying equation
      \cite{Car97,BR03}, without a-priori knowledge of how the final mesh will look like.
      In a similar vein,
      for large-scale computations with more than a few tens of millions of
      elements and massively parallel systems, the I/O related to creating and accessing the
      mesh data structure is often a serious bottleneck. Much better performance
      can be obtained by reading smaller meshes that get refined
      as part of the simulation.
      In all of these cases, refining the mesh involves the computation of new grid
      points from inside the simulation, and these points need to respect
      the same geometry used for the original mesh.
    \item Accurate simulators use curved cells and higher order mappings both at
      the boundary and in the interior of the mesh. How exactly these curved cells
      should look requires information about the underlying
      geometry. To illustrate the importance of this point, let us
      mention that it has been understood theoretically for a long time that one loses
      the optimal order of finite element discretizations with
      polynomial degrees greater than one, if one does not also use
      higher-order approximations of the boundary
      \cite{Ciarlet1972,Mansfield1978,Bassi1997,Bra97,Bartels2004,Dor98,BN99}.
      This is also known from practical experience; for example
      \cite{Har02} presents experimental evidence, and
      \cite{mengaldo2020industryrelevant} and the references therein
      provide an excellent example of the lengths one needs to go to
      to recover higher-order accuracy if the underlying geometry is
      not available.
      (Table~\ref{tab:spherical} and Fig.~\ref{fig:transfinite} below also
      illustrate the issue.) In other words, a finite element solver that does not know
      about the underlying geometry will compute a solution
      with an unnecessarily large error or use an unnecessarily large
      amount of computational work. Alternatively, additional points for a high-order representation
      need to be computed as a separate workflow step between the meshing and the
      simulation, as used, e.g., in \cite{Hindenlang2015,Krais19,Moxey2020}.
    \item In many contexts -- during the simulation, but also during accurate
      visualization and evaluation of quantities of interest -- it is important
      to know the correct normal vector to faces of cells. An approximation can
      be obtained by simply taking the location of vertices as provided in the
      mesh file and their connection into cells, as the ground truth. But the
      vectors so computed are not consistent with the true, underlying geometry
      from which the mesh was originally generated. As a consequence, visualizations
      are often not faithful representations of the actual computations, and
      quantities of interest are not evaluated to the full accuracy
      possible. For example, accurate evaluation of mass or energy fluxes across a boundary requires accurate knowledge of the normal vector.
      For fourth order equations, accurate knowledge of the
      normal vector may also be required to retain optimal convergence
      order of finite element schemes if it is necessary to construct
      $C^1$ approximations of a curved boundary (see
      \cite{Mansfield1978} and the references therein).
\end{itemize}
We will provide more examples below where geometry information is used in finite
element simulations.

These considerations point to a need to propagate geometry information not only to the mesh generator, but
indeed also to the simulation engine. This raises the question how this
can best be done. To the best of our knowledge, no commercial or open source tools
do this in a consistent way today. Furthermore, we have to realize that geometries are
often described through complex CAD systems that are either not open source, complicated
to interact with, or can only be accessed in proprietary ways; as a
consequence, it makes sense to ask what kinds of queries a generic geometry
engine has to be able to answer to satisfy the needs of simulation software. In order to address all
of these points, we have undertaken a comprehensive study with the following
goals:
\begin{enumerate}
    \item Identify a \textit{minimal set of operations} -- which we will call ``primitives'' -- 
      that geometry tools need to be
      able to perform to satisfy the needs of simulation tools.
    \item Provide a \textit{comprehensive review of geometry operations} performed
      in a widely used finite element library and a large application code
      that is built on the former, with the aim of verifying that the
      minimal set of operations outlined above is indeed sufficient.
    \item Discuss \textit{ways in which geometry tools can implement the minimal
      set of operations}, given the kinds of geometries and information that
      is typically available in industrial and research workflows.
\end{enumerate}
 We will state the primitives in the form of \textit{oracles}, i.e., as blackboxes that given
 certain inputs produce appropriate outputs, without specifying how exactly
 one would need to implement this operation. This reductionist approach is 
 often appropriate when one wants to interface with one of many possible
 geometry engines, each of which may have its own way of implementing the operations.
 In such cases, it is often useful to only specify a minimal interface that
 all engines can relatively easily implement. A common way of representing
 oracles in object-oriented codes is by providing an abstract base class with
 unimplemented virtual functions; the base class is then the oracle, whereas
 derived classes provide actual implementations.

The result of our work is the realization that only two geometry primitives
are sufficient. We will discuss these in
Section~\ref{sec:primitives} and will empirically show in Section~\ref{sec:uses}
that they can indeed satisfy (almost) all needs of simulators. (The sole exception
requiring additional information from the geometry engine is discussed in
Appendix~\ref{sec:exceptions}). Section~\ref{sec:implementing-the-primitives} is then
devoted to the question of how one would implement these two operations in
the most common situations, namely where the geometry may be described explicitly
(for example, if the geometry is a sphere or a known perturbation of it) or where
it is described implicitly (e.g., through a collection of NURBS patches in
typical CAD engines). We demonstrate the practical benefits of the integration of 
geometry and simulation in Section~\ref{sec:examples}. We briefly discuss
periodic domains in Appendix~\ref{sec:periodicity}.

The practical implication of our work is the identification of a minimal
interface that allows the coupling of geometry and simulation engines.
We have tested these interfaces via the widely used finite element software \dealii{} \cite{dealiicanonical,dealII91}
and the Advanced Simulator for Problems in Earth ConvecTion \aspect{}~\cite{KHB12,HDGB17}
with geometry descriptions that are either given explicitly, or
via the OpenCASCADE library that is widely used for CAD descriptions \cite{opencascade}.
All of the results of our work are available in the publicly released
versions of \dealii{} and \aspect{}.
The verification of our approach using these examples, and the fact
that the sufficient interface is so small, should provide the certainty
necessary to follow a similar path for integrating other simulation and geometry
software packages.

We end this introduction by remarking that one may also wish to provide
geometry information to the postprocessing stage. For example, this would
allow visualization software to produce more faithful graphical representations
of the solutions generated by simulators, free of artifacts that result from
incomplete knowledge of the domain on which the simulation was performed.
We are not experts in visualization and consequently leave an investigation
of the geometry needs of visualization software to others. We will, however,
comment that the evaluation of quantities of interest -- such as stresses
at individual points, heat and mass fluxes across boundaries, or average and
extremal values for certain solution fields -- may be most efficiently
done from within the simulator itself, given that geometry information as
well as knowledge of shape functions and other details of the discretization
are already available there; indeed, this is the approach chosen in the
\aspect{} code mentioned above.

\section{Fundamental geometric primitives}
\label{sec:primitives}

As we will discuss in detail in the following section, it turns out
that the geometric queries needed by finite element codes
-- such as finding locations for new vertices upon mesh refinement, or
computing vectors normal to the surface -- can be reduced to only two
``primitive'' operations: (i) finding a new point given a set of existing
points with corresponding weights,
and (ii) computing the tangent vector to a line connecting
two existing points. This realization of a minimal set of
operations allows us great flexibility in choosing software packages
that actually provide these primitives, and minimizes the dependency
a finite element code incurs when using an external geometry package.

In the literature, implementations that can answer a small number of
very specific questions -- i.e., provide certain simple operations --
are typically referred to as ``oracles''. The 
point of postulating the existence of an oracle is that it allows us
to separate the \textit{design} of a code from its actual
\textit{implementation}. In the current case, all that matters for the
purposes of the current section is that an oracle exists that can
answer two specific questions, and whose answers can be used
throughout a finite element code.

In order to motivate the two geometric primitives that we postulate are sufficient
for almost all finite element operations, let us first provide two
scenarios of relevance to us. First, consider a
$d$-dimensional surface embedded into a higher dimensional space; one might think
of this surface as the boundary of a volume within
which we would like to simulate certain physics. The second setting
is a $d$-dimensional volume geometry in $d$-dimensional space
for which we would like to consider interior cells to also deviate
from the simplest, $d$-linear shape; an example would be a volume mesh
that extends away from a curved wing around which we would like to simulate
air flow. We will use these two scenarios for the examples below.


\subsection{Statement of primitives}

Given this background, the two operations we have found are necessary
are the following:

\begin{primitive}[``New point'']
  Given $N\ge 2$ existing points $\mathbf x_1,\ldots,\mathbf x_N$ and weights
  $w_1,\ldots,w_N$ that satisfy $\sum_{n=1}^N w_n=1$, return a new
  point $\mathbf x^\ast(\mathbf x_1,\ldots,\mathbf
  x_N,w_1,\ldots,w_N)$ that interpolates the given
  points, weighting each $\mathbf x_n$ with $w_n$.
\end{primitive}

\begin{primitive}[``Tangent vector'']
  Given two existing points $\mathbf x_1,\mathbf x_2$, return the
  (non-normalized) ``tangent'' vector $\mathbf t$ at point $\mathbf x_1$ in direction
  $\mathbf x_2$, defined by 
  \begin{align}
    \label{eq:tangent}
    \mathbf t(\mathbf x_1,\mathbf x_2) 
    = \lim_{w\rightarrow 0}
    \frac{\mathbf x^\ast(\mathbf x_1,\mathbf x_2, (1-w), w)-\mathbf x_1}{w}.
  \end{align}
\end{primitive}

We will give many examples below of how these two operations are
used. As we will show, all other geometric questions a typical finite
element code may have can be answered exactly, without approximation, with only these two operations.
For the moment, one can think of use cases as follows:
\begin{itemize}
    \item When adaptively refining a mesh, one needs to introduce new
    vertex locations on edges, faces, and in the interior of cells. This
    is easily done using the first of the two primitives above, using the
    existing vertex locations on an edge, face, or cell as input.
    
    \item Computing the normal vector to a face at the boundary of a three-dimensional
    domain can be done by taking the cross product of two tangent vectors that pass
    through the point at which we need the normal vector. These tangent
    vectors are computed via the second primitive.
\end{itemize}
A concrete implementation that describes a particular geometry will be able to
provide answers to the two operations based on its knowledge of the actual
geometry. For example, and as discussed in
Section~\ref{sec:implementing-the-primitives}, the wing of an airplane would be described
using a CAD geometry that can be queried for the primitives above at least
on the boundary of the domain. More work -- also described in Section~\ref{sec:transfinite} below -- will
then be necessary to extend this description into the interior. In other cases,
however, the description of the geometry may be available everywhere -- for example, 
if the geometry of interest is an analytically known object such as a sphere.

\begin{remark}
  \label{remark:no-domain}
  As we will see below, the operations that finite element codes require
  can be implemented by only querying information from objects that
  describe an \textit{entire} manifold, without having to know
  anything about the triangulation that lives on it, or in particular
  where the triangulation's boundary lies. This greatly simplifies the
  construction of oracles because they do not need to know anything
  about meshes, or the domain on which we solve an equation. For
  example, we will be able to use a CAD geometry of, say, the
  entire hull of a ship even if we want to solve an equation on only
  parts of the surface; that there is a boundary to the domain on
  which we solve, and where on the hull it lies, will be of no
  importance to the oracle that will answer our queries.

  Information about the domain, its boundary, and the mesh that covers
  it, will of course be used in constructing the
  \textit{inputs} to the queries, but is not necessary in computing their
  \textit{outputs}.
\end{remark}

\begin{remark}
  \label{remark:tangent-finite-differences}
  Given the definition of the tangent vector in \eqref{eq:tangent}, it
  is possible to implement this second query approximately
  through finite differencing with the help of the first, using
  \begin{align*}
    \mathbf t 
    \approx
    \frac{\mathbf x^\ast(\mathbf x_1,\mathbf x_2, (1-\varepsilon), \varepsilon)-\mathbf x_1}{\varepsilon}
  \end{align*}
  with a small but finite $\varepsilon$. In other words, one could
  in principle get away with only one primitive if necessary.
  At the same time, and as
  discussed in Section~\ref{sec:implementing-the-primitives},
  we have found that in many cases, good ways to directly implement
  the \textsc{tangent vector} primitive exist, and there is no need
  for approximation.
\end{remark}

\begin{remark}
  \label{remark:minimality}
  The two primitives mentioned above are the only ones necessary to
  implement the operations discussed in the next section, exactly and
  without approximation. As such, they are truly \textit{primitive},
  but that does not mean that one could not come up with a larger set
  of operations geometry packages could provide, possibly in a more
  efficient way than when implemented based on the primitives. We did
  not pursue this idea any further, primarily because (i)~we have not
  found it necessary in practice, and (ii)~because these additional
  operations would have to be implemented for all of the different
  ways discussed in
  Section~\ref{sec:implementing-the-primitives}. (That said, see also Remark~\ref{remark:third-primitive} in the appendix.)
\end{remark}


\section{A survey of geometry uses in finite element codes}
\label{sec:uses}

In our initial attempt to understand what kinds of primitives suffice to implement typical
finite element operations, and later to support our claim that the two primitives
mentioned above are indeed sufficient for almost all operations,
we have undertaken a survey of a large and modern finite element
library to find all of the places in which finite element codes make
use of geometric information. Specifically, we have investigated the
\dealii{} library (see \url{https://www.dealii.org}, \cite{dealiicanonical,dealII91}), consisting of more than
a million lines of C++. \dealii{} provides many modern
finite element algorithms, for example fully adaptive meshes,
higher order elements and mappings, support for coupled and nonlinear
problems, complex boundary conditions, and many other areas. We have also surveyed the more than 150,000 lines
of C++ code of the \aspect{} code that simulates convection in the
Earth mantle \cite{KHB12,HDGB17}.

We have found a large number of places in which these codes use
geometric information, but most of them can be grouped into general
categories. In the following subsections, we present a summary of use cases
along with a description of how each of these uses of
geometric information can be implemented in terms of only the two
primitive operations outlined above. We have found only one exception
of an operation -- rarely used in finite element operations -- that
can not be mapped to the two primitives and would need additional
information from the geometry engine that underlies the computation;
we discuss it in Appendix~\ref{sec:exceptions}.

The intent of this summary is to show the capabilities enabled by
these two primitives. Given their simplicity, we expect them to
be deliverable by any software package that
provides information about the domain on which a partial differential
or integral equation is to be solved. This includes, in particular,
the many CAD packages that are used today to model essentially every
object around us.

\subsection{Mesh refinement}

When refining a mesh, either adaptively or globally, one typically
splits each edge into two children, and then proceeds with faces and
cell interiors in ways that depend on whether we use
triangles/tetrahedra or quadrilaterals/hexahedra; both cases follow
essentially equivalent schemes.

\begin{figure}
  \includegraphics[width=0.23\textwidth]{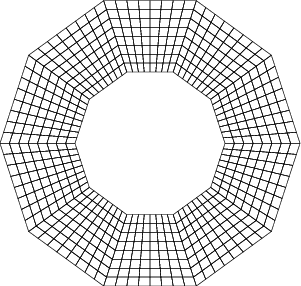}
  \hfill
  \includegraphics[width=0.23\textwidth]{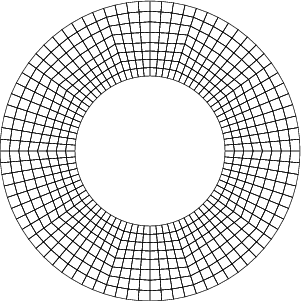}
    \hfill
  \includegraphics[width=0.23\textwidth]{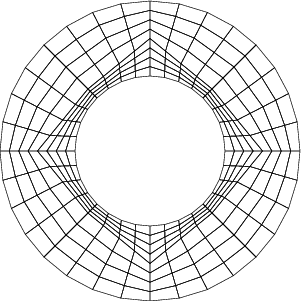}
  \hfill
  \includegraphics[width=0.23\textwidth]{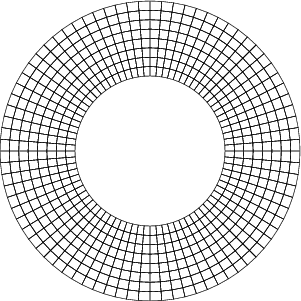}
  \caption{\it Illustration of the importance of taking geometry into
    account when refining meshes. Starting from a discretization of an annulus
    that contains ten coarse cells, we choose new points either ignoring the geometry description, i.e., computing the location of a new vertex by simply averaging
    the locations of the surrounding points (left), or choosing  new points on boundary
    edges so that they have the correct radius from the origin, and
    new points of interior edges and cells as the Cartesian mean of the surrounding vertices (second from left).  The latter procedure works well when the number of coarse elements is sufficient to resolve the geometry, but leaves some kinks in the grids, where one is still able to identify the original ten coarse cells. However, it may lead to very distorted grids if the coarse mesh is not fine enough, e.g., if the coarse mesh consisted of only four cells (second from right). The ideal case (right) is independent of the number of coarse cells that one may start with, and it exploits full knowledge of the underlying geometry. In this case, this is done by choosing all new points on edges and cells so that
    they average the \textit{radius and angle} of the adjacent points.}
  \label{fig:hypershell}
\end{figure}

For breaking an edge, the new midpoint needs to
respect the geometry of the underlying domain. For example, if the
existing edge is at the boundary of a domain, then we will want the
new point to also lie on the boundary; its location should
then not simply be at the mean 
of the Cartesian locations of the existing vertices of the edge. The same
is often true for \textit{interior} edges: if refining a mesh for a
domain around an obstacle, e.g., a wing, then one typically also wants
to refine interior edges so that they use properties of the
geometry -- in other words, we want the \textit{entire} mesh to
follow the geometry of the wing, not just the faces or cells
actually adjacent to the wing. 
We illustrate this using a simpler example: refining a mesh in a spherical shell
geometry (see Fig.~\ref{fig:hypershell}). There, 
we will want to introduce new points on interior edges so
that they have \textit{radius} and \textit{angle} equal to the average of the existing
points. A similar strategy must be used for the new point
at the center of each triangular or quadrilateral facet, whether at the
boundary or interior, and then again for the new point inside each
tetrahedral or hexahedral cell if in 3d; in these cases, the interpolation
happens between more than just two adjacent vertex locations.

In all of these cases, the new point on edges, faces, and cells can
clearly be satisfied with the \textsc{new point} primitive of
Section~\ref{sec:primitives}. (The different
meshes shown in Fig.~\ref{fig:hypershell} then simply correspond
to different ways of implementing the primitive.) In the case of isotropic
refinement of edges, one needs to call the primitive with the two
existing points and equal weights $w_1=w_2=\frac 12$. For the new
point at the center of a quadrilateral, we first compute
the new points on the edges, and then call the \textsc{new point}
primitive for the new center point with the 8 surrounding points (the
four pre-existing vertices of the quadrilateral plus the four new
points on the now split edges). Assigning the existing vertices a weight
of~$-1/4$ and the new edge centers a weight~$1/2$ leads to high-quality
meshes, as opposed to putting equal weights on all points. These
weights correspond to the evaluation with
a transfinite interpolation \cite{Gordon82} in the center of the edges.
Similarly, for hexahedra, we call the \textsc{new point} primitive
with the 8 existing vertices, the 12 new edge points, and the 6 new
face centers, with weights~$1/8$, $-1/4$, and~$1/2$,
respectively. Similar considerations can also be applied to triangles
and tetrahedra.

\subsection{Polynomial mappings}
\label{sec:mappings}

Many finite element codes use ``iso-parametric'' mappings in which
cells consist of the area of the reference cell (e.g., the unit
square) mapped by a polynomial mapping of degree equal to that of the
finite element in use. (One may of course also use lower, ``sub-parametric'', or even higher order,
``super-parametric'', mappings.)

\begin{figure}
    \centering
    \includegraphics[width=0.7\textwidth]{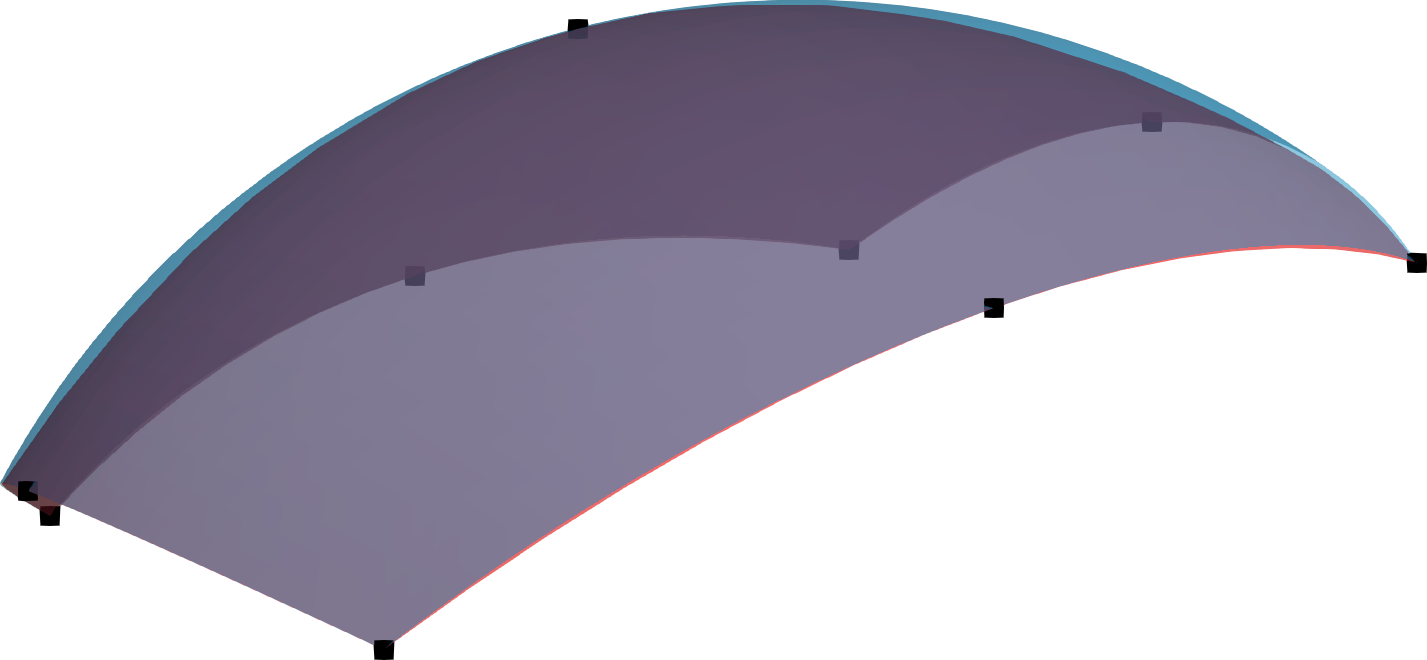}
    \caption{\it Illustration of biquadratic mapping of a quadrilateral (shown in red) onto the surface of a sixth of a sphere (blue), with interpolation nodes represented as black squares.}
    \label{fig:quadratic_mapping}
\end{figure}

When using quadratic or higher order mappings, cells are bounded by polynomial
curves or surfaces. In general, for quadrilaterals and hexahedra, a
cell $K$ is defined by
\begin{align}
  \label{eq:polynomial-mappings}
  K &= \left\{
    \mathbf x \in {\mathbb R}^d : \mathbf x 
    = \mathbf F_K(\mathbf {\hat x}) = \sum_i \mathbf v_i
    \varphi_i(\mathbf {\hat x}), \quad \mathbf {\hat x} \in \hat K=[0,1]^d\right\}.
\end{align}
Here, $\varphi_i$ are the usual $d$-linear, $d$-quadratic, or higher
order shape functions, and $\mathbf v_i$ are the $(p+1)^d$ support
points of cell $K$ when working with mappings of degree $p$. (Similar constructions apply to triangles and tetrahedra.)

For $p=1$, i.e., the bi-/trilinear mappings, the support points are
simply the vertices of the cell, for which we know that they are
already part of the manifold, whose locations are fixed, and which are
generally provided as part of the mesh description. On the
other hand, for higher order mappings, we need to evaluate the
locations of these support points between vertices. For example, for
cubic mappings with equidistant support points, one needs these points
at relative distances of 1/3 and 2/3 along each edge. This is easily
done using the \textsc{new point} primitive, using the two vertices at
the end of the edge as existing points, and weights $w_1=2/3, w_2=1/3$
for the first support point, and reversed weights for the
second. Following the computation of edge support points, we proceed
with the evaluation of support point locations inside each face
using the already computed points at the vertices and on
the edges as pre-existing points; finally, we compute locations for
points inside cells. The weights within the entities might again be derived
from transfinite interpolation \cite{Gordon82}, propagating the
information from surrounding edges into the interior.
Figure~\ref{fig:quadratic_mapping} illustrates the construction of such an
approximate, polynomial representation of a part
of a sphere. The support points for the mapping that are not existing
vertices
are generated using the \textsc{new point} primitive; consequently, they lie exactly
on the curved (spherical) manifold, whereas there is a (small) gap between the
polynomial approximation and the original sphere everywhere else.

All operations outlined in the previous paragraph only require the
first of our two primitives.
In Section~\ref{sec:uses-c1}, we analyze the requirements for implementing
$C^1$-mappings that are needed in certain finite element applications. 

\subsection{Computing the Jacobian of a mapping}
\label{sec:jacobian}

The computation of integrals in the finite element method generally involves the transformation
$\mathbf F_K: \hat K \mapsto K$ from the reference cell $\hat{K}$ to the concrete cell
$K$. The transformation of the integrand then implies that we need the
``Jacobian'' $J_K=\hat\nabla \mathbf F_K(\mathbf {\hat x})$ or its inverse for all vector quantities, and the
determinant, $\text{det}(J_K)$, as an additional weight factor. The
Jacobian matrix $J_K$ is also necessary when using a Newton iteration
to find the reference coordinates $\mathbf{\hat x}=\mathbf F^{-1}_K(\mathbf x)$ that correspond
to a given point $x$; this operation is key in particle-in-cell
methods and methods based on characteristics.
In the
following, let us assume that $\mathbf F_K$ maps to the ``true'' domain
of interest, not a polynomial approximation of it; we will comment on the
latter case at the end of the section.

The construction of the matrix $J_K(\mathbf {\hat x})$ at a specific quadrature
point $\mathbf {\hat x}_0$ given in reference coordinates uses an
algorithm that first builds a basis for the tangent space at $\mathbf x_0=\mathbf F_K(\mathbf {\hat x}_0)$. There are of course many bases for
this tangent space, but we choose the one whose basis vectors are the
images of unit vectors $\mathbf{\hat e}_i$ centered at $\mathbf{\hat x}_0$
as the images of these vectors will correspond to the partial derivatives
$\frac{\partial}{\partial \hat x_i} \mathbf F(\mathbf{\hat x}_0)$ and
consequently form the columns of $J$. The algorithm proceeds in three steps:
\begin{enumerate}
\item Use the \textsc{new point} primitive with arguments given by the vertices of the concrete cell, and weights equal to the reference coordinates of $\mathbf{\hat x}_0$, to obtain the image $\mathbf{x}_0 \in K$ of the point $\mathbf{\hat x}_0 \in \hat K$;
\item Select $d$ points $\{\mathbf{\hat x}_i\}_{i=1}^d$ in the reference cell $\hat K$, aligned with the $i$th reference coordinate axis, and compute their images in $K$ using again the \textsc{new point} primitive with arguments given by the vertices of the cell $K$, and weights equal to the coordinates of the points $\mathbf{\hat x}_i$;
\item Call the \textsc{tangent vector} primitive $d$ times, with arguments $(\mathbf x_0, \mathbf x_i)$, divide the result by $L_i = \mathbf{e}_i\cdot(\mathbf{\hat x}_i-\mathbf{\hat x}_0) = (\mathbf {\hat x}_i-\mathbf{\hat x}_0)_i$,
to obtain tangent vectors $\mathbf t_i$.
\end{enumerate}
These vectors $\mathbf t_i$ form a basis of the tangent space (though they may
neither be orthogonal nor normalized) and are then arranged to form the columns of
$J_K(\mathbf{\hat x}_0)$.

\begin{figure}
	\centering
	\includegraphics[]{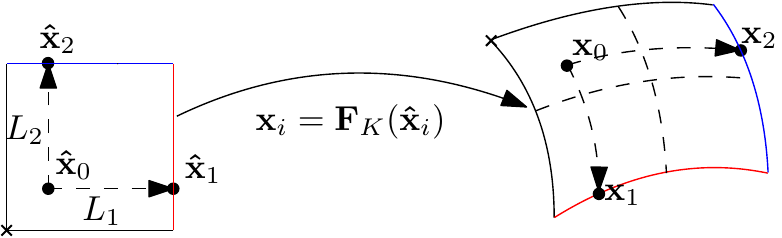}
\caption{\it Computing the Jacobian of a mapping (at point $\mathbf{\hat x}_0$) using only the \textsc{New Point} and \textsc{tangent vector} primitives.}
\label{fig:compute-jacobians}
\end{figure}

While one might in principle choose the points $\mathbf{\hat x}_i$
arbitrary in either $\pm\mathbf{\hat e}_i$ from $\mathbf{\hat x}_0$, in practice one wants to choose them as far
away from $\mathbf{\hat x}_0$ as possible in the reference cell
to obtain a well conditioned algorithm (see Fig. \ref{fig:compute-jacobians}).

We note that in many cases, finite element codes do not use the \textit{exact}
mapping $\mathbf F_K$ from the reference cell $\hat K$ to a real cell $K$,
but only a polynomial approximation (see Section~\ref{sec:mappings}). In those
cases, it is of course trivial to compute the Jacobian matrix since it
is simply the gradient of the mapping presented in \eqref{eq:polynomial-mappings},
i.e., 
\begin{align*}
    J_K(\mathbf {\hat x}) = \sum_i \mathbf v_i
    \hat\nabla\varphi_i(\mathbf {\hat x}).
\end{align*}
It is an exercise to show that the two methods indeed result in the same matrix for polynomial mappings.

\subsection{Computing vectors normal to lower-dimensional manifolds}
\label{sec:normal-vectors}

In our survey of existing uses of geometric information, we have found
many instances where normal vectors to lower-dimensional surfaces
embedded in higher dimensional space were required. These could either
be the normals to $(d-1)$-dimensional faces of $d$-dimensional cells
in $d$ space dimensions, or the normals to $d$-dimensional cells
embedded in a $d'$-dimensional space where $d'>d$.
Examples of such situations include the following:
\begin{itemize}
\item \textit{In the evaluation of boundary conditions:} Across
  applications, there exist a wide variety of boundary conditions that
  require normal vectors. Among these are:
  \begin{itemize}
    \item ``No flux'' boundary conditions on a vector velocity field
      $\mathbf u$ have the form $\mathbf n \cdot \mathbf u = 0$. These
      require knowledge of the normal vector in the implementation because one
      typically imposes a constraint that determines one vector
      component in terms of the others, for example by evaluating the
      identity $u_3 = -\frac{n_1 u_1 + n_2 u_2}{n_3}$ at the support
      points of the shape functions.
    \item In other cases, one may want to prescribe a particular
      tangential velocity in addition to no-flux conditions such as
      the one above. An example is the classical lid-driven cavity. In
      this case, one needs to also constrain the tangential part of
      the velocity, $(\mathbf I-\mathbf n \otimes \mathbf n)\mathbf u$.
    \item For hyperbolic equations, boundary conditions can only be
      imposed on inflow boundaries. Whether a part of the boundary is
      in- or outflow depends on the sign of $\mathbf n \cdot \mathbf
      u$ where $\mathbf u$ is either an externally prescribed
      velocity, or the solution of a previous time step or nonlinear
      iteration. Again, an explicit construction of the normal vector
      to the boundary is required.
    \item In fluid-structure interaction problems,
      traction boundary conditions on the solid are common where the
      traction $T$ is given by $\mathbf T=p\mathbf n$ in terms of the
      fluid pressure $p$. It may also contain a viscous friction
      component $2\eta \varepsilon(\mathbf u)\mathbf n$ where $\eta$
      is the viscosity, $\mathbf u$ is the fluid velocity and
      $\varepsilon(\mathbf u)$ its symmetric gradient.
    \item In external scattering problems, one often only solves for
      the scattered, instead of the total field. The boundary
      conditions on the scatterer then involve the incident field
      $U_\infty$, and depending on the type of scatterer often contain
      terms of the form $U_{\infty}\mathbf n$.
  \end{itemize}

\item \textit{In discontinuous Galerkin (DG) discretizations:} When using
  the discontinuous Galerkin method, one generally applies integration by
  parts on each element
  individually, resulting in face integrals of some quantities times the normal
  vector on the face. The numerical fluxes that weakly impose the continuity of
  solutions between the elements also often involve the normal vector~\cite{Hesthaven07}. 
  One example are upwind formulations for hyperbolic equations, which 
  consider the sign and magnitude of the quantity $\mathbf n\cdot\mathbf u$
  where again $\mathbf u$ is a velocity field. Thus,
  the normal vector explicitly appears
  in the bilinear form that defines the discretization.

\item \textit{In error estimation:} Many error estimators for finite
  element discretizations contain terms that measure the ``jump''
  from one cell to the next of the normal component of the gradient of
  the discrete solution. An example is the widely used ``Kelly'' error
  estimator for the Laplace equation \cite{KGZB83,GKZB83} that
  requires evaluating
  \begin{align*}
    \eta_K = \left(\frac{h}{24} \int_{\partial K} \left[\mathbf n
      \cdot \nabla u_h \right]^2 \; ds \right)^{1/2}
  \end{align*}
  for every cell $K$, where $[\cdot]$ denotes the jump across a cell
  interface. Similar terms appear in almost all other residual-based
  error estimates (see, for example, \cite{AO00,BS01book,BR03}).

\item \textit{In boundary element methods:} Boundary integral formulations hinge upon the availability of a fundamental solution for the PDE at hand; i.e., an (often singular) function $G$, such that the solution $u$ of a PDE at a point $\mathbf x\not\in\Gamma$ can be written
as a function of the jump of the solution and of its normal derivative on the boundary $\Gamma$:
\begin{equation*}
    u(\mathbf{x}) + \int_\Gamma \mathbf{n}_y \cdot \nabla_y G(\mathbf{x}-\mathbf{y})  [u(\mathbf{y})]d\Gamma_y =  \int_\Gamma [\mathbf{n}_y \cdot \nabla_y u(\mathbf{y})] G(\mathbf{x}-\mathbf{y}) d\Gamma_y;
\end{equation*}
see, for example,~\cite{Sauter2011,Mola2013,GiulianiMolaHeltai-2015-a,Mola2017}.

\item \textit{In postprocessing solutions:} Computer simulations are
  typically done because we want to \textit{learn} something from the
  solution, i.e., we need to evaluate or postprocess it after it is
  computed. Many of these postprocessing steps involve normal vectors
  at the boundary. For example, in \aspect{} alone, we evaluate the normal
  components of a flow field (or of the stress tensor) at the boundary
  to displace the boundary appropriately; we may compute the gravity
  field and its angle against actual surface topography to determine
  the direction of water and sediment transport; and we compute the
  normal component of the gradient of the temperature field to
  determine the heat flux through surfaces. Many other postprocessing
  applications come to mind that require the normal vector to the
  boundary. In other contexts, we require the normal vector to an
  embedded surface on which we solve equations, for example to assess
  the evolution of surfaces such as soap films or the membranes of
  cells.
\end{itemize}

In all of the cases above -- typical of many complex finite
element simulators -- we require computing the normal vector to faces
or cells. The construction of the normal vector $\mathbf n\in {\mathbb R}^d$
follows in essence the
construction outlined in Section~\ref{sec:jacobian} for the Jacobian matrix, except that we map a lower-dimensional object (namely, a face of a cell) and
can consequently only generate $d-1$ basis vectors. The algorithm then finds
a set of $d-1$ tangent vectors $\mathbf t_i$ to the manifold (this involves
using the \textsc{tangent vector} primitive), and we compute the normal
vector using the wedge product,
\begin{align}
  \label{eq:normal}
    \mathbf n = \pm \frac{\mathbf t_1 \wedge \cdots \wedge \mathbf t_{d-1}}
                         {\|\mathbf t_1 \wedge \cdots \wedge \mathbf t_{d-1}\|}.
\end{align}
The sign is typically chosen so that $\mathbf n$ points \textit{outward} from
whatever object we are currently considering. The algorithm in
Section~\ref{sec:jacobian} chooses a particular basis, that is, in
particular neither orthogonal nor normalized. Neither is a problem here
because the definition of $\mathbf n$ above is invariant under the choice
of basis.

There are two situations that require additional thoughts:
\begin{itemize}
\item If one asks for the normals of faces of cells that are
  themselves embedded in a higher dimensional space: For example, a
  user code may require the jump of the normal derivative across the
  one-dimensional edges of two-dimensional cells, where the cells form
  a triangulation of a two-dimensional manifold embedded in
  three-dimensional space. In such a case, two orthogonal normals
  exist to the edge, and the one required in the code will likely be
  the one that is tangent to the manifold.

  In general, objects of dimensionality $d_1\le d_M$ that are part of a
  triangulation of a manifold $M$ of dimension $d_M$ embedded in $d\ge d_M$
  dimensional space allow for a system of $d-d_1$ mutually orthogonal
  normal vectors. In applications, we typically seek those that also
  lie in the $d_M$-dimensional tangent space to the manifold, leaving us
  with a choice of $d-d_1-(d-d_M)$ mutually orthogonal normal vectors.

\item In a similar case, we frequently seek the normal vector to a
  triangulated manifold of dimension $d_M$ embedded in a $d>d_M$
  dimensional space. An application is to determine the direction of
  growth for the two-dimensional membrane surrounding a
  three-dimensional biological cell. In other applications, however, we consider
  embeddings in even higher dimensional spaces, $d>d_M+1$, in which
  case there is again a system of $d-d_M$ mutually orthogonal normal
  vectors that we will need to be able to compute.
\end{itemize}
In both of these situations, we can only compute a subset of the complete
set of tangent vectors $\mathbf t_i, i=1,\ldots,d-1$; the normal vector
appropriate for the application can then generally be computed by providing
the remaining tangent vectors from the context to complete the wedge product
in \eqref{eq:normal}.

\subsection{Computing vectors tangent to lower-dimensional manifolds}

There are also cases where computing tangent vectors to a boundary is
necessary. For example, in electromagnetics applications, one often requires
that the tangent component of the electric field vanishes. Mathematically,
this is often expressed by requiring that 
$\mathbf n(\mathbf x) \times \mathbf E(\mathbf x)=0$ for all points
$\mathbf x$ on the boundary of a domain, but
actual applications generally implement this by requiring that
$\mathbf t_1(\mathbf x) \cdot \mathbf E(\mathbf x) 
  = \cdots 
  = \mathbf t_{d-1}(\mathbf x) \cdot \mathbf E(\mathbf x) = 0$,
where the vectors $\mathbf t_i(\mathbf x)$ form a complete (but otherwise arbitrary)
basis of the tangent space to the boundary at the point $\mathbf x$. Each of these
conditions is then imposed in the same way as the constraints corresponding
to $\mathbf n \cdot \mathbf u = 0$ in the previous section.

Since the concrete choice of basis vectors $\mathbf t_i$ is unimportant, we can again fall back to the construction discussed in
Section~\ref{sec:jacobian} in order to impose 
$\mathbf t_i(\mathbf x) \cdot \mathbf E(\mathbf x) = 0$. 
On the other hand, for stability purposes, it may also be convenient
to orthonormalize these vectors, which is of course easily done using,
for example, the Gram--Schmidt process.

\subsection{$C^1$ mappings}
\label{sec:uses-c1}

As mentioned above, one frequently uses polynomial mappings from the reference
cell $\hat K$ to each concrete cell $K$ of the mesh. In practice, this
implies that curved boundaries are approximated by piecewise polynomials. The common construction of these mappings guarantees the continuity
of the individual pieces of this boundary approximation at vertices, but the
approximate boundary is not $C^1$ continuous at vertices (in 2d) and edges
(in 3d). On the other hand, in some applications, it is desirable or necessary
that curved boundaries are approximated with $C^1$ continuity. 

One possible way
to achieve this is to use cubic mappings, and to construct their face support points such that the mapping interpolates the actual manifold at the face vertices in such a way that the resulting mapping's tangent space in vertices coincides with the tangent space of the underlying, original, curved boundary.
If the boundary itself is continuously differentiable across face boundaries, this will guarantee that the resulting mapping is also continuously differentiable. To achieve this, we need to require that the cubic
mapping $\mathbf F_K$ on cell $K$ satisfied the following conditions in each vertex $\mathbf v_i, i=1,\ldots,2^{d-1}$ of the face:
\begin{align}
        \mathbf{F}_K(\mathbf{\hat v}_i) = \mathbf{v}_i, \qquad &&
        J_K(\mathbf{\hat v}_i)\mathbf{\hat n}_i = \mathbf{n}_i,
    \label{eq:C1-mapping-continuity-constraints}
\end{align}
where the normals $\mathbf{n}_i$ are computed as in Section~\ref{sec:normal-vectors}, while the normals $\mathbf{\hat n}_i$ are unit vectors perpendicular to the faces of the reference cell. 

The above conditions can then be used to fix the locations of the (non-vertex) support points of the cubic mapping on the face, by providing $2^{d}$ conditions on the $2^d$ support points on the faces. (However, these non-vertex locations
will, in general, not actually lie on the original boundary.) The remaining support points can be set using the \textsc{new point} primitive, as in Section~\ref{sec:mappings}, completing the information necessary for the
construction of a cubic mapping $\mathbf F_K$.

\section{Geometry specification formats}
\label{sec:geometry-formats}

In order to understand how the two primitives can be implemented in common, industrial use cases, let us use this section to discuss the ways in which
geometry information is typically available: namely, either in the form of
combinations of simple geometries, or in the form of a collection of analytically known
patches that together provide a boundary representation
for the domain of interest.

\subsection{Constructive Solid Geometry}
\label{sec:csg}

In some cases, the shape of the computational domain admits an
analytic representation for which it is straightforward to implement
the primitives described in Section~\ref{sec:primitives}. For example,
if the domain
corresponds to a cube, a parallelogram, or a sphere, then it can be described
in the language of differential geometry where each local portion of the
domain (a ``chart'') is the image of a Cartesian space under an analytically
known ``push-forward'' operation. In such cases, the \textsc{new point}
operation is easily represented by applying the ``pull-back'' to the
surrounding points, averaging in the Cartesian domain, and ``pushing forward''
the resulting point. The \textsc{tangent vector} operation follows
the same idea. We discuss these algorithms in more detail
in Section~\ref{sec:chart-case}.

It is straightforward to generalize this
procedure for shapes that are further transformed from relatively simple
original geometries -- for example, when using a true topographic model of
the Earth instead of a sphere; in those cases, the pull-back and push-forward
operations are simply concatenations of the Cartesian-to-simple-geometry
and simple-to-transformed-geometry operations.

\begin{figure}
  \phantom{.}\hfill
\includegraphics[width=.35\textwidth]{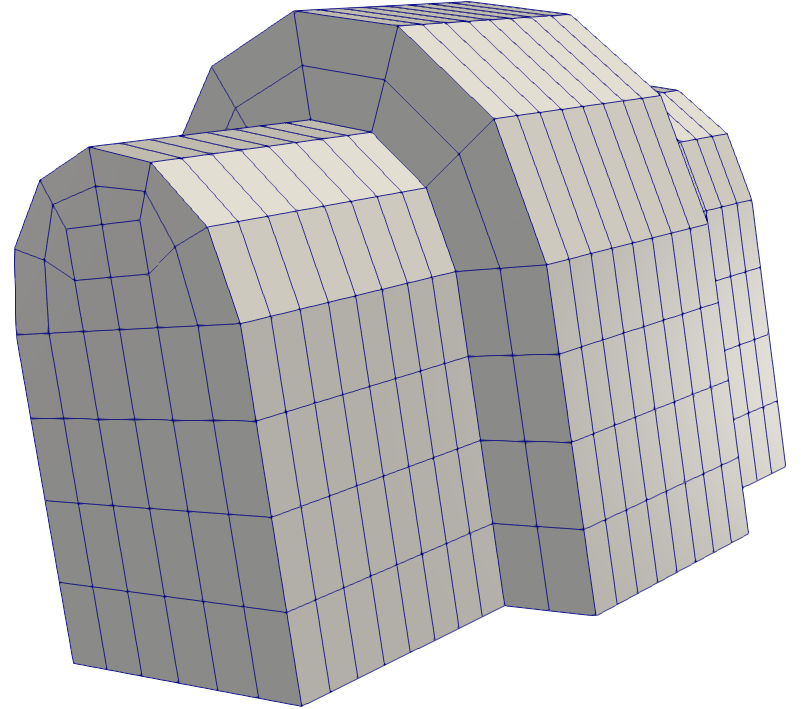}\hfill
\includegraphics[width=.35\textwidth]{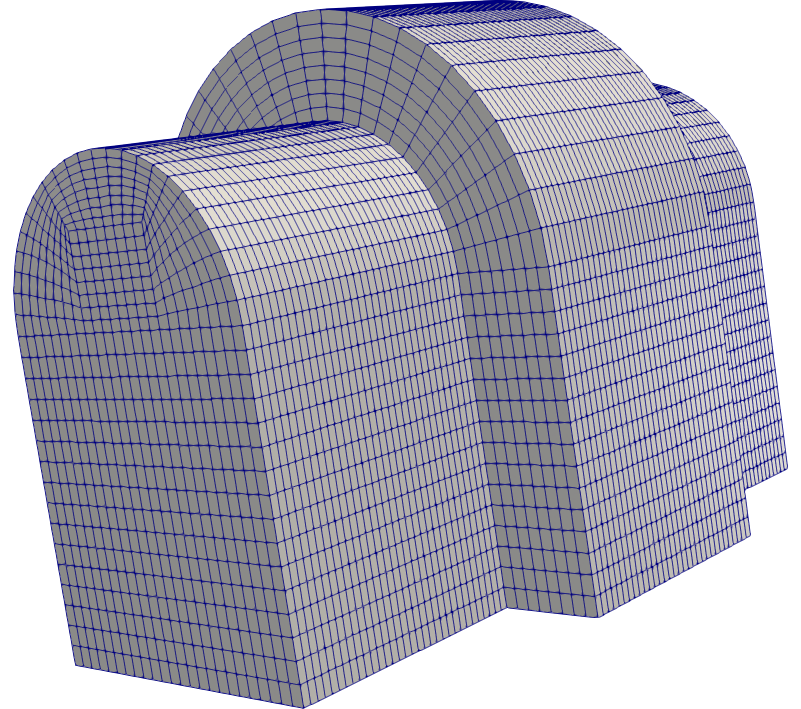}
  \hfill\phantom{.}
\caption{\it Example of geometry constructed using simple Constructive
  Solid Geometry primitives and shapes,
  taken from the \textsc{deal.II} tutorial program \texttt{step-49} \cite{step-49}.
  Left: Coarse mesh. Right: A refined mesh
that uses the underlying geometry description of the domain's pieces to
place new vertices.
}
\label{fig:step-49}
\end{figure}

A generalization of this approach is commonly known as ``Constructive Solid Geometry'' (CSG): Here, solids are expressed as as the union, intersection,
or set difference of simpler solids (e.g., cubes, prisms, spheres, cylinders, cones, and toruses, possibly subject to transformations).
An example involving boxes and cylinders is shown in Figure~\ref{fig:step-49}.
The implementation
of our primitives
then requires keeping track of which part of the geometry one is on, and
the primitives are implemented by falling back to the case of a single,
analytically known geometry for which we have provided a description
above.%
\footnote{In practice this generally requires that the mesh respects
the boundaries between the original elementary geometry building blocks
so that the primitives are only ever called with points belonging to the
same geometric part. Most mesh generation programs do not provide
such a guarantee unless provided surfaces separating pieces of the geometry, for example if the combined
geometry is a \textit{non-overlapping} union of elementary geometrical shapes.}

\subsection{Boundary Representation Models}
\label{sec:brep}

In the vast majority of industrial cases, 
domains are specified through a boundary representation model (B-REP or BREP), used 
to combine topological and geometrical information to provide a
complete description of the boundary of the domain. A solid (or surface) is then represented as a collection of connected surface (or curve) elements:
this set of
$d-1$-dimensional surface(s), embedded into $d$-dimensional space, represents the
boundaries of the domain.

There are two general
ways in which these boundaries are often provided: (i) In CAD applications,
curves and surfaces are generally
parameterized as Non-Uniform Rational B-Splines (NURBS); (ii)
The surface of the domain is given as
a triangulated surface, often obtained from a point cloud
derived from measurements on an actual object.%
\footnote{In actual practice, both the NURBS ``patches'' and the triangles produced by widely
  used software often overlap, contain gaps, or leave holes, leading to
  difficult questions of interpretation that we will ignore here. See
  also \cite{mengaldo2020industryrelevant}.}

\subsubsection{Boundaries represented through NURBS}
Here,
individual patches are defined as parametric surface patches that
have an analytical representation. In three space dimensions, this representation is
often provided by a mapping from a ``reference domain'' with coordinates $u,v$
to a NURBS surface. So, within each patch, one might think that 
all information associated with the primitive operation is easily obtained using this
mapping from $u$-$v$ space to ${\mathbb R}^3$. In practice, however, this mapping may map
equal areas in $u$-$v$ space to rather unequal areas on the surface in ${\mathbb R}^3$;
one also encounters difficult issues if one needs to cross from one NURBS patch to another
when evaluating our two primitives. We will therefore discuss the practical implementation of the
primitives in Section~\ref{sec:implementing-the-primitives} below, and for the moment
simply state that the primitives \textit{can} be implemented based on
NURBS B-REPs.

Practical realizations need to read and
properly interpret the content of CAD data structures. This is made possible through open source
libraries such as OpenCASCADE~\cite{opencascade}, which provide tools for importing and interrogating vendor-neutral
CAD file formats such as IGES and STEP \cite{IGES-5.3,Marjudi2010,STEP}.

\subsubsection{Boundaries represented by triangulated point clouds}
\label{sec:point-cloud-formats}
In many other cases, the boundary is described by points:
a polygonal chain in 2d, or a triangulation in 3d.
These are often used to represent the most topologically intricate domains, and are for instance
used by software which extracts geometrical models from MRI images or 3d laser scans. Several CAD
modeling tools can also export triangulated surfaces. A popular file
format is STL, which originated in the
stereolithography community but is now widely used in industrial continuum mechanics applications and
3D printer geometry specifications.

Despite their wide usage, triangulated surfaces also pose several problems in the implementation
of suitable primitive operators. In principle, the primitives can locate new points
and compute tangent vectors based on the planar surfaces associated
with each triangle. However, there is typically no guarantee that the
triangles form a closed surface -- models may have
holes or gaps in which a new point or a tangent vector cannot
be identified. More importantly, in order to provide effective information to a finite element algorithm, the diameters of triangles must be significantly smaller than those of the finite element mesh:
Only then will the triangulated surface appear
smooth compared to the finite elements mesh. Such requirements lead in most
situations to surfaces triangulated by extremely large numbers of elements, resulting
in high storage costs and expensive computations to identify new points and tangent vectors
as both require looping over the many triangular faces.

\section{Possible implementations}
\label{sec:implementing-the-primitives}

Section~\ref{sec:uses} has shown that essentially all geometrical queries required in finite element codes can be implemented by only using the two operations 
described in Section~\ref{sec:primitives}. However, we have not discussed
how these primitives can be implemented for the kinds of
geometries one frequently encounters in simulations. Section~\ref{sec:geometry-formats}
might have given some hints, but also outlined some of the difficulties.
This section discusses common approaches.

It is important to note that the way this is done is not unique, and there are often different
possible implementations of the primitives. For example, in Section~\ref{sec:examples} we will
show results where $\mathbf x^\ast$ is chosen as the directional,
orthogonal, or closest point projection onto a surface. These choices may
result in meshes of different quality, but this is immaterial to the
key point of this paper as they only affect how an oracle
\textit{implements} a query, not \textit{the kinds of queries} it has to be
able to answer.

\subsection{Implementing primitives for analytically known charts}
\label{sec:chart-case}
The simplest way to implement the two primitives is if we have direct and explicit access
to an atlas that describes the actual manifold of which the domain is a part. 
We can then rely on the language and methods of differential geometry.

An atlas for a manifold $M$ of dimension $d$ is a collection of coordinate ``charts'' (also called coordinate patches, or local frames): homeomorphisms $\phi_\alpha $ from an open subset $U_\alpha$ of $M$ to an open subset $U_\alpha^\ast$ of the Euclidean space of dimension $d$, such that $\cup_{\alpha} U_\alpha = M$.
By definition, the sets $U_\alpha$ and $U_\alpha^\ast$ are connected by push-forward and
pull-back functions, $\phi_{\alpha}: U_\alpha^\ast \rightarrow U_\alpha \subseteq M$ and
$\phi_{\alpha}^{-1}: U_\alpha \rightarrow U_\alpha^\ast$.

With these definitions, if a point $\mathbf x_2$ lies in the same set $U_\alpha$ of the point $\mathbf x_1$ (i.e., if
they lie on the same \textit{chart}), then the parametric line that
connects them is the image of the straight line
$\phi_\alpha^{-1}(\mathbf x_1)+t\left(\phi_\alpha^{-1}(\mathbf x_2)-\phi_\alpha^{-1}(\mathbf x_1)\right)$ in $U_\alpha^\ast$ under
$\phi_\alpha$, i.e., 
\begin{align*}
  \mathbf s(t)
  =
  \phi_\alpha\left(\mathbf x_1^\ast+t\left[\mathbf x_2^\ast-\mathbf x_1^\ast\right]\right)
  =
  \phi_\alpha\left(\phi_\alpha^{-1}(\mathbf
  x_1)+t\left[\phi_\alpha^{-1}(\mathbf x_2)-\phi_\alpha^{-1}(\mathbf x_1\right]\right).
\end{align*}

Therefore, when the push-forward and the pull-back are explicitly available, a straightforward implementation of the  \textsc{new point} primitive with any number of points that belong to the same set $U_\alpha$ is given by the following
    \begin{align*}
  \mathbf x^\ast(\mathbf x_1,\ldots,\mathbf
  x_N,w_1,\ldots,w_N)
  =
  \phi_\alpha\left(\sum_{n=1}^N w_n \phi_\alpha^{-1}(\mathbf
  x_n)\right).
\end{align*}
That is, we simply take the weighted average of the
pulled back points $\phi_{\alpha}^{-1}(\mathbf x_n)$ and push the
average forward.

Furthermore, the \textit{tangent vector}
primitive is also easy to compute using the formula
\begin{align*}
  \mathbf t(\mathbf x_1,\mathbf x_2)
  =
  \nabla^\ast\phi_\alpha\left(\mathbf x_1^\ast\right)
  \left[\mathbf x_2^\ast-\mathbf x_1^\ast\right]
  =
  \nabla^\ast\phi_\alpha\left(\phi_\alpha^{-1}(\mathbf
  x_1)\right)
  \left[\phi_\alpha^{-1}(\mathbf x_2)-\phi_\alpha^{-1}(\mathbf x_1)\right],
\end{align*}
where $\nabla^\ast \phi_\alpha(\mathbf x^\ast)$ is the gradient of the push-forward
function for the chart around $\mathbf x^\ast$.

\subsection{Implementing primitives based on geodesics}
\label{sec:geodesic-case}

The formulas available for analytically known charts suggest that the quality of the implementation depends on whether or not $\phi_\alpha$ is constructed in a ``reasonable'' way.
In order to produce the least distorted grids upon refinement, one would like $\phi_\alpha$ to map straight lines in $U_\alpha^\ast$ to geodesic paths on $M$. In general, the construction of explicit geodesics is only possible on simple manifolds, and even there, subtle problems lurk.

For example, it is tempting to construct the \textsc{new point} primitive with two points directly in terms of \textit{geodesics}. If we parameterize a geodesic that connects $\mathbf x_1, \mathbf
x_2$ as $\mathbf s(t)$ so that $\mathbf x_1=\mathbf s(0)$, $\mathbf
x_2=\mathbf s(1)$, and assume that $\mathbf s(t)$ moves at constant
speed (as measured by the metric), then this suggests the following implementation: given $\mathbf
  x_1,\mathbf x_2$ and weights $w_1,w_2=1-w_1$, define
  \begin{align*}
    \mathbf x^\ast\left(\mathbf x_1,\mathbf x_2,w_1,1-w_1\right)
    = \mathbf s(w_1),
  \end{align*}
Similarly, we can construct the \textsc{tangent vector} primitive for two points $\mathbf x_1,\mathbf x_2$ connected by the geodesic $\mathbf s(t)$ as
  \begin{align*}
    \mathbf t
    = \mathbf s'(0).
  \end{align*}

While reasonable, this construction solely based on geodesics is only useful if geodesics connecting two points are unique. In general, this requires that the points are ``close together''
in some sense -- which in the finite element context means that the mesh is already sufficiently
fine since we generally call the primitives with points that are located on the same cell.
Moreover, its generalization to more than two points is not at all trivial.
Indeed, consider a \textit{recursive} implementation for the \textsc{new point} primitive with more than two points: Given
  $\mathbf x_1,\ldots,\mathbf x_N$, $w_1,\ldots,w_N$, $N>2$, define
  \begin{equation*}
    \mathbf x^{\ast}\left(\mathbf x_1,\mathbf x_2,\ldots,\mathbf x_N,w_1,w_2,\ldots,w_N\right)
    = 
    \mathbf x^\ast\left(\mathbf
    s_{12}\left(\frac{w_1}{w_1+w_2}\right),\mathbf x_3\ldots,\mathbf x_N,
           w_1+w_2,w_3,\ldots,w_N\right)
  \end{equation*}
  where, without loss of generality, we have assumed that point
  $\mathbf x_1$ has a nonzero weight $w_1>0$, and
  where $\mathbf s_{12}(t)$ is the (presumed unique) geodesic connecting $\mathbf x_1,\mathbf
  x_2$. It is not difficult to show that,
  in general, the algorithm above depends on the order in which points and weights are given.
  Furthermore, the operation is not associative in the following sense:
 Consider, for example, a situation with the four vertices of a quadrilateral cell and four equal weights $w_i=\tfrac 14$. Using the recursive definition of $\mathbf x^\ast$ above, one can compute
\begin{align*}
        \mathbf x_{1234} & =  \mathbf x^{\ast}(\mathbf x_1, \mathbf x_2, \mathbf x_3, \mathbf x_4, \tfrac 14, \tfrac 14, \tfrac 14, \tfrac 14) \\
        \tilde{\mathbf x}_{1234} & =  \mathbf x^{\ast}(\mathbf x_{12}, \mathbf x_{34}, \tfrac 12, \tfrac 12)
        \quad\text{with }
        \mathbf x_{12}  =  \mathbf x^{\ast}(\mathbf x_1, \mathbf x_2, \tfrac 12, \tfrac 12)
        \text{ and }
        \mathbf x_{34} =  \mathbf x^{\ast}(\mathbf x_3, \mathbf x_4, \tfrac 12,\tfrac 12).
\end{align*}
In general, $\tilde{\mathbf x}_{1234} \neq 
\mathbf x_{1234}$, contrary to expectation. However,
both are ``reasonable'' intermediate points, and it is not
clear which option to prefer. This ambiguity suggests that this may not be
a useful algorithm.

A commutative alternative algorithm
is to take the average over all possible permutations of the pairs $(\mathbf x_i, w_i)$:
\begin{equation*}
  \mathbf x^{\ast}\left(\mathbf x_1,\mathbf x_2,\ldots,\mathbf x_N,w_1,w_2,\ldots,w_N\right) := \frac{1}{N!}\sum_{k = 1}^{N!}\mathbf x^{\ast}_\text{recursive}(S_k(\mathbf x, w)),
 \end{equation*}
where $S_k(\mathbf x, w)$ is the $k$th permutation of the $N$ pairs $\{\mathbf x_i, w_i\}_{i=1}^N$
and $\mathbf x^\ast_\text{recursive}$ is the recursive implementation from before.
Of course, such an algorithm is, in general, quite expensive given that, for example, one
already has to consider $8!=40,320$ permutations for finding the point that interpolates 
the vertices of a hexahedron.

There is also the issue of how exactly one finds the appropriate
point on a geodesic. In many practical applications, a geodesic is easy
to describe geometrically, but parameterizing it is more complicated;
in those cases, one may be tempted to first average the input points in
the ambient space and only then \textit{project onto the geodesic} -- avoiding to
provide the geodesic with an arc-length that can be used to find intermediate
points. But this yields yet another complication: Consider, for example,
\begin{equation}\label{eq:refine_identity}
        \mathbf{x}_{12} = \mathbf{x}^{\ast}(\mathbf x_1, \mathbf x_2,
        \tfrac 34, \tfrac 14),
        \qquad\qquad\qquad
        \tilde{\mathbf{x}}_{12} = \mathbf{x}^{\ast}\big(\mathbf x_1, \mathbf{x}^{\ast}(\mathbf x_1, \mathbf x_2, \tfrac 12, \tfrac 12), \tfrac 12, \tfrac 12\big).
\end{equation}
The first formula represents the point
at $1/4$ of the geodesic between $\mathbf x_1$ and $\mathbf x_2$, whereas the second first computes
the mid point $\mathbf{x}_\text{m}$ between $\mathbf{x}_1$ and $\mathbf {x}_2$, and then the mid point between $\mathbf{x}_1$ and $\mathbf{x}_\text{m}$. Unfortunately, the definition of $\mathbf x^\ast$
based on the projection onto geodesics does not always provide that the two points are the same.%
\footnote{An example is the unit circle embedded into $\mathbb R^2$, where
$\mathbf x^\ast$ first averages the input points in the ambient space with
their weights, and only then projects back onto the circle. In this
context, consider $\mathbf x_1=(1,0), \mathbf x_2=(0,1)$. Then 
$\mathbf{x}^{\ast}(\mathbf x_1, \mathbf x_2, \tfrac 34, \tfrac 14)$
first computes the weighted average of the input points, yielding $(\tfrac 34,\tfrac 14)$, which is then projected onto the circle: 
$\mathbf x_{12}\approx (0.9487,0.3162)$. On the other hand,
$\mathbf x_m=\mathbf{x}^{\ast}(\mathbf x_1, \mathbf x_2, \tfrac 12, \tfrac 12)
 = (0.7071,0.7071)$ and
$\tilde{\mathbf{x}}_{12}\approx(0.9239,0.3827)$ -- a different point.
A similar situation of course happens on the surface of a sphere embedded
in 3d.}
It is not
difficult to construct examples where this ambiguity has concrete, detrimental effects on the
accuracy of finite element operations. For example, when interpolating a finite element
solution from a parent cell $K$ (say, the image of the interval $(0,1)$ under a transformation
$\Phi_K$) to its first child $K^{\text{child 0}}$ (the image of $(0,\tfrac 12)$ under $\Phi_K$),
then we will want that the point with reference coordinates $\xi$ on the child cell equals
the point with reference coordinates $\xi/2$ on the parent cell when evaluating finite element
shape functions. This can, in general, be ensured for $d=1$, but not for $d\ge 2$: There, the
construction above (where we have considered $\xi=\tfrac 12$)
shows that that the two points may be different, and one can observe that this affects the 
convergence order of some algorithms.


\begin{table}
    \centering
    \caption{\it Effect on the numerical error for interpolation between two successive mesh levels in a spherical shell, using a spherical manifold with spherical averages as the \textsc{new point} oracle.}
    \label{tab:spherical}
    \begin{tabular}{lllclll}
    \hline
    \multicolumn{3}{l}{Polynomial degree $p=4$} & & \multicolumn{3}{l}{Polynomial degree $p=7$}\\
    \hline
    $n_{\text{DoF,coarse}}$ & error coarse & error after refine & & $n_{\text{DoF,coarse}}$ & error coarse & error after refine \\
    490 & $7.25\cdot 10^{-2}$ & $6.92\cdot 10^{-2}$ && 2,368 & $2.92\cdot 10^{-3}$ & $9.70\cdot 10^{-3}$ \\
    3,474 & $2.42\cdot 10^{-3}$ & $4.16\cdot 10^{-3}$ && 17,670 & $1.16\cdot 10^{-5}$ & $3.28\cdot 10^{-3}$ \\
    26,146 & $1.05\cdot 10^{-4}$ & $2.74\cdot 10^{-4}$ && 136,474 & $6.68\cdot 10^{-8}$ & $2.56\cdot 10^{-4}$ \\
    202,818 & $3.38\cdot 10^{-6}$ & $1.66\cdot 10^{-5}$ && 1,072,626 & $2.21\cdot 10^{-10}$ & $1.62\cdot 10^{-5}$ \\
    1,597,570 & $1.06\cdot 10^{-7}$ & $1.02\cdot 10^{-6}$ && 8,505,058 & $8.71\cdot 10^{-13}$ & $1.02\cdot 10^{-6}$ \\
    \hline
    \end{tabular}
\end{table}

For the sphere, where the construction of geodesics is trivial, \cite{Buss01} provides exhaustive discussions that show that 
various possibilities for evaluating $\mathbf x^\ast$ with three or more input points 
yield results that
differ by quantities of order ${\mathcal O}(D^4)$ where $D$ is a measure of the distance between the 
input points. Table~\ref{tab:spherical} shows the result of an experiment in this
vein: Given a
subdivision of a spherical shell in ${\mathbb R}^3$ with inner and outer radii $0.5$ and $1.0$ into hexahedra, we define a
finite element field $u_h$ of polynomial degree $p$ whose nodal values are chosen so that
$u_h$ interpolates a known smooth function $u$. We show $\|u-u_h\|_{L_2}$ in the columns labeled ``error coarse''
for $p=4$ and $p=7$ in the table. We can observe the error decrease approximately 
as ${\mathcal O}(h^{p+1})$ as the mesh is chosen finer and finer, as expected.

For each of these meshes, we then perform one isotropic refinement and interpolate $u_h$ onto
$u_{h/2}$ by using the embedding of the two spaces \textit{on the reference cell}
(rather than interpolating $u$ onto $u_{h/2}$). Of course, one would expect the
difference $\|u-u_{h/2}\|$ to be the same as before since we expect that $u_{h/2}=u_h$.
But this equality assumes that the support
points of the child cells of size $h/2$ (when computed using the interpolation between
the vertices of the child cell) are located at the positions where one would have expected
them to be when interpolating the (child cell's) support points using the vertices of
the parent cell. Given considerations such as in \eqref{eq:refine_identity}, this is
not the case, and due to this difference, the error $\|u-u_{h/2}\|$ is actually substantially
\textit{larger} than $\|u-u_h\|$ and, instead of the optimal convergence rate $\mathcal O(h^{p+1})$,
only converges as $\mathcal O(h^4)$.%
\footnote{It will not surprise readers that one can spend many hours on ``debugging'' a code 
  as the identity $u_{h/2}=u_h$ is clearly ``obvious'' given how $u_{h/2}$ was constructed.
  One can restore optimal convergence by directly
  interpolating $u_h$ onto $u_{h/2}$ using the node points on the
  actual child cells of the sphere, instead of pulling $u_h$ back to
  the reference cell for each cell, computing the interpolation to its children, and then
  pushing forward from the child cell in reference coordinates to the child cells of
  the refined mesh.}

These, and the examples below, demonstrate that the implementation of $\mathbf x^\ast$
is not obvious even in situations where we have reasonable pull-back and
push-forward operations. Indeed, the message of this section is that
there are no \textit{obvious} solutions, but a variety of approaches that
yield \textit{reasonable} implementations good enough for many cases -- the
fact that they might affect convergence rates when used without a
deeper understanding notwithstanding.

There are also cases where not even the definition of
pull-back and push-forward functions is
easily possible. For example, if the points  $\mathbf x_n$ do not fall within the \textit{same} set $U_\alpha$, it may not be possible to construct a \textit{single} chart parameterized by $\phi_{\alpha}$ that
encompasses \textit{all} points $\mathbf x_n$. This may lead to unpleasant
ambiguities, such as in the case of periodic domains (see
Appendix~\ref{sec:periodicity}) or where only piece-wise descriptions of the geometry are available
(e.g., CAD geometries via collections of NURBS patches).

\subsection{Implementing primitives based on projections onto CAD geometries}
\label{sec:projection-based}

For CAD geometries, any two points may fall \textit{across} two different (non-overlapping) NURBS
patches $U_\alpha$ and $U_\beta$. But there are more difficulties that make it difficult to
implement the primitives based on pull back/push forward approaches: 
\begin{enumerate}
\item NURBS patches do not always satisfy the requirements of a chart
  in the topological sense, i.e., they may be non-invertible in some
  points, and they may be non-smooth ( i.e., contain corners and edges
  \textit{within} a single $U_\alpha$).
\item The metric $\phi_\alpha$ that results from mapping $u$-$v$ space to a NURBS patch $U_\alpha$
may be ill-formed: Its derivative may be close to zero near some
points, and large at others --
equally spaced points in $u$-$v$ space (or points $\mathbf x^\ast(\mathbf x_1, \mathbf x_2; w_1, (1-w_1))$ for equally spaced values of $w_1$) would then be mapped to highly unevenly spaced points.
\item The evaluation of the pull-back $\phi_\alpha^{-1}(\mathbf x)$ is computationally very expensive for NURBS patches.
\end{enumerate}
As a consequence, \textit{useful} implementations of our primitives will not be based on the
concepts of differential geometry, but will rather be \textit{projection-based}.
Indeed, a possible implementation of the \textsc{new point} primitive for (multi-patch) CAD-based
geometries is provided by the following algorithm:
Given $\mathbf x_1,\ldots,\mathbf x_N$, $w_1,\ldots,w_N$, $N\ge 2$, define
  \begin{equation*}
    \mathbf x^{\ast}\left(\mathbf x_1,\mathbf x_2,\ldots,\mathbf x_N,w_1,w_2,\ldots,w_N\right)
    = 
    P_{\text{CAD}}\left( \sum_{n=1}^{N} w_n \mathbf x_n\right),
    \end{equation*}
where $P_{\text{CAD}}(\mathbf x)$ is a projection of the point $\mathbf x$ onto the CAD surface (or
curve). We will discuss various choices for the projection operator
$P_{\text{CAD}}$ below,
given that the choice will 
influence the quality of the result. 
We note that projection-based strategies have been proposed in \cite{Farin02} and are used
in a basic variant for a high-order finite element code in \cite{Krais19}.
Many CAD programs, and specifically the OpenCASCADE library~\cite{opencascade} that we use for the examples shown here,
implement all of the operations necessary for the three approaches
to implementing a projection discussed below.

\subsubsection{Projection in a fixed direction}
\label{sec:direction-based}
The cheapest way to compute a projection is if the direction of the
projection $\mathbf d\in {\mathbb R}^d$ is known a priori -- for example, because we know that
the surface in question has only minor variation from being horizontal.
In that case, one might choose
  \begin{equation*}
    \mathbf x^{\ast}\left(\mathbf x_1,\mathbf x_2,\ldots,\mathbf x_N,w_1,w_2,\ldots,w_N\right)
    = 
    I_{\text{CAD}}\left( \sum_{n=1}^{N} w_n \mathbf x_n, \mathbf d \right),
    \end{equation*}
where $I_{\text{CAD}}(\mathbf x, \mathbf d)$ is the intersection of the line
$\mathbf s(t)=\mathbf x + s\mathbf d$ and the CAD surface (or curve).
In other words, we first average all points $\mathbf x_i$ and then 
move the average back onto the surface along direction $\mathbf d$.

\subsubsection{Projections taking into account the normal vector}
\label{sec:normal-to-mesh-based}
In more general situations, however, choosing a projection direction
a priori is not possible. Rather, one needs to take into account the
geometry of the CAD surface in the vicinity of the point to be projected,
for example by considering the normal vectors to either the
existing mesh, or to the surface.

\begin{figure}
  \centering
  \includegraphics[width=.5\textwidth]{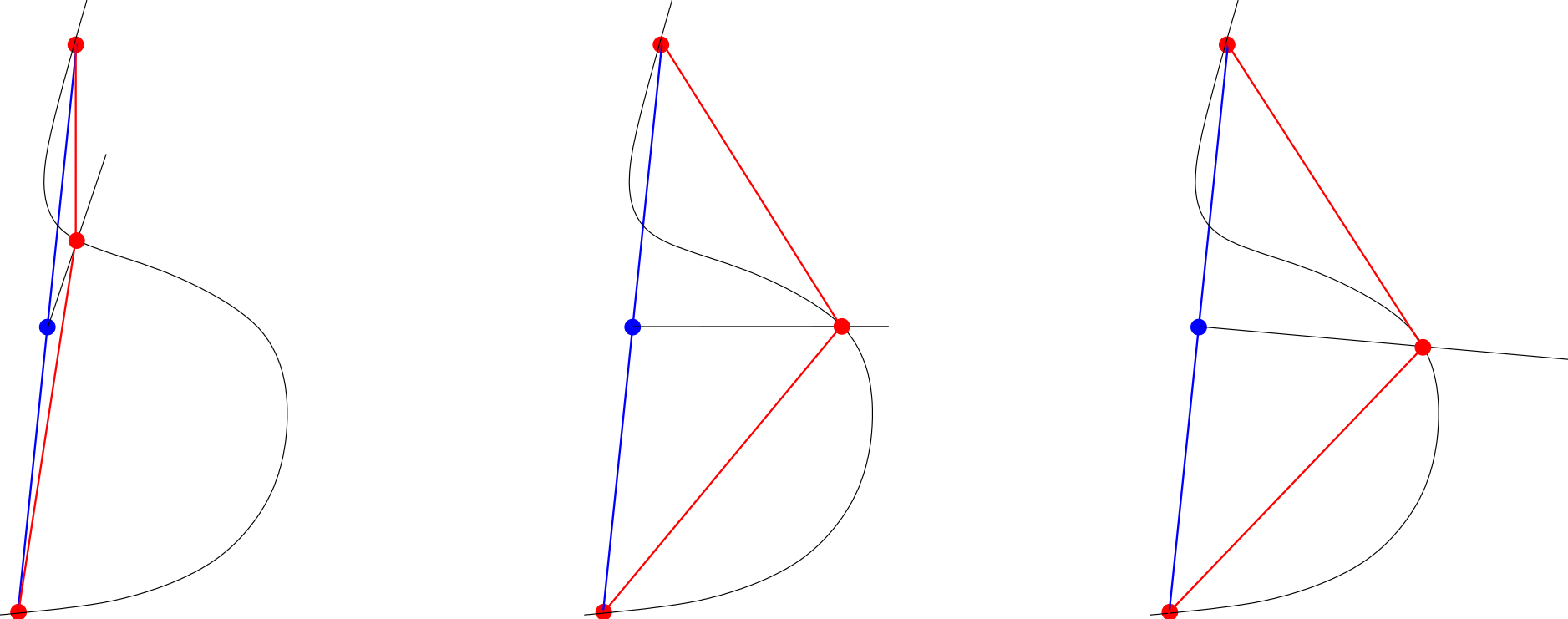}
  \caption{\it Comparison of three different implementations of the \textsc{new point} primitive for CAD geometries. The red end points (coarse vertices)
  of the blue line (coarse cell) form the inputs $\mathbf x_1,\mathbf x_2$
  for which we want to use the \textsc{new point} primitive to find
  a new mid-point (i.e., $w_1=w_2=\tfrac 12$). The blue point is the average of the original vertices to be projected onto the curved
  geometry. Left: Projection normal to the geometry. Center: Projection in
  a direction chosen a priori. Right: Projection normal to coarse mesh.}
  \label{fig:surface-projection-comparison}
\end{figure}

The first of these options (shown in the right panel of \ref{fig:surface-projection-comparison}) would use the following
implementation:
  \begin{equation*}
    \mathbf x^{\ast}\left(\mathbf x_1,\mathbf x_2,\ldots,\mathbf x_N,w_1,w_2,\ldots,w_N\right)
    = 
    I_{\text{CAD}}\left( \sum_{n=1}^{N} w_n \mathbf x_n, \mathbf n \right),
    \end{equation*}
where the direction $\mathbf n$ is now (an approximation) to the normal
vector of the area identified by the points $\mathbf x_n$. $\mathbf n$
is clearly defined if one only has $d$ input points in $d$ dimensional
space; if there are more points -- e.g., the four vertices of a face of
a hexahedron in 3d -- then one will want to define some useful approximate
vector, e.g., the vector normal to the least squares plane that approximates
the point locations. As is clear from the figure, if this
implementation is chosen for mesh refinement, one generally ends up with
refined meshes with cells of rather uniform sizes.

To use this approach, one needs to have at least $d$ input points
in $d$ dimensions in order to define a unique direction normal to
the existing points. But we also need to be able to use the
\textsc{new point} primitive when finding a new midpoint for an edge
in 3d. Thus, we need an additional condition to identify a unique
direction among all of those perpendicular to the line connecting the
existing edge end points. We do this by averaging the CAD surface normal at the vertices of the edge (both of which we know are on the
surface), and then projecting it onto the edge's axial plane.

An alternative, often implemented in CAD tools but expensive to evaluate,
is to use a direction vector $\mathbf n$ that is \textit{perpendicular to
the actual geometry}, rather than to the current mesh. This is
shown in the left panel of Fig.~\ref{fig:surface-projection-comparison}
and may lead to child cells of different size. Ultimately, however,
once the mesh is already a good approximation to a surface, both
of the approaches mentioned here will yield very similar results.

\subsubsection{Implementing the \textsc{tangent vector} primitive for
  CAD surfaces}

In all three of the project-based cases above, the implementation of
the \textsc{tangent vector} primitive may be constructed using a
finite differences approximation as already mentioned in
Remark~\ref{remark:tangent-finite-differences}. A more accurate
approach pushes forward the tangent vector at the pulled-back
point.

\subsection{Extending boundary representations into volumes}
\label{sec:transfinite}

The previous sections only dealt with finding
points and tangent vectors on a lower-dimensional surface, given
points already on that surface. On the other hand, finite element codes
typically use volume meshes for which the CAD geometry only provides
information about the \textit{boundary}. Thus, one still needs
a way to extend this information into the \textit{interior} of the
domain -- the importance of this step is apparent by looking at 
Fig.~\ref{fig:hypershell}.

A general mechanism for this task is based on transfinite
interpolation~\cite{Gordon82}. A transfinite interpolation maps points from some
reference space $\hat{\mathbf x}\in \hat K=[0,1]^d$ to points in real space $\mathbf x$
by a weighted sum of information on the geometry of the faces of the image of $\hat
K$. For example, for a quadrilateral
in two dimensions
\begin{multline*}
    \mathbf{x}(\hat x_1, \hat x_2) = (1-\hat x_2)\mathbf c_0(\hat x_1)+\hat x_2 \mathbf c_1(\hat x_1) +
    (1-\hat x_1)\mathbf c_2(\hat x_2) + \hat x_1 \mathbf c_3(\hat x_2) \\
    - \left[(1-\hat x_1)(1-\hat x_2) \mathbf x_0 + \hat x_1(1-\hat x_2) \mathbf x_1 + (1-\hat x_1)\hat x_2 \mathbf x_2 + \hat x_1 \hat x_2 \mathbf x_3 \right].
\end{multline*}
Here, $\bf c_0(s), \bf c_1(s), \bf c_2(s), \bf c_3(s)$ are the
four parameterized curves describing the geometry of the edges of the deformed
quadrilateral and $\bf x_0, \bf x_1, \bf x_2, \bf x_3$ are the four vertices. The
evaluation on each edge is
done via the \textsc{new point} primitive $\mathbf{x}^{\ast}$, i.e.,
$\mathbf c_0(s) = \mathbf{x}^{\ast} (\mathbf x_0, \mathbf x_1, 1-s, s)$
and similarly for the other curves. If an edge is straight, then $\bf c_0(s) =
(1-s)\mathbf{x}_0 + s \mathbf{x}_1$. Similar formulas extend to
the three-dimensional case. The important point is that transfinite mappings
exactly respect the geometry of the boundary, while extending it smoothly
into the interior of the domain.

\begin{figure}
    \centering
    \includegraphics[width=0.22\textwidth]{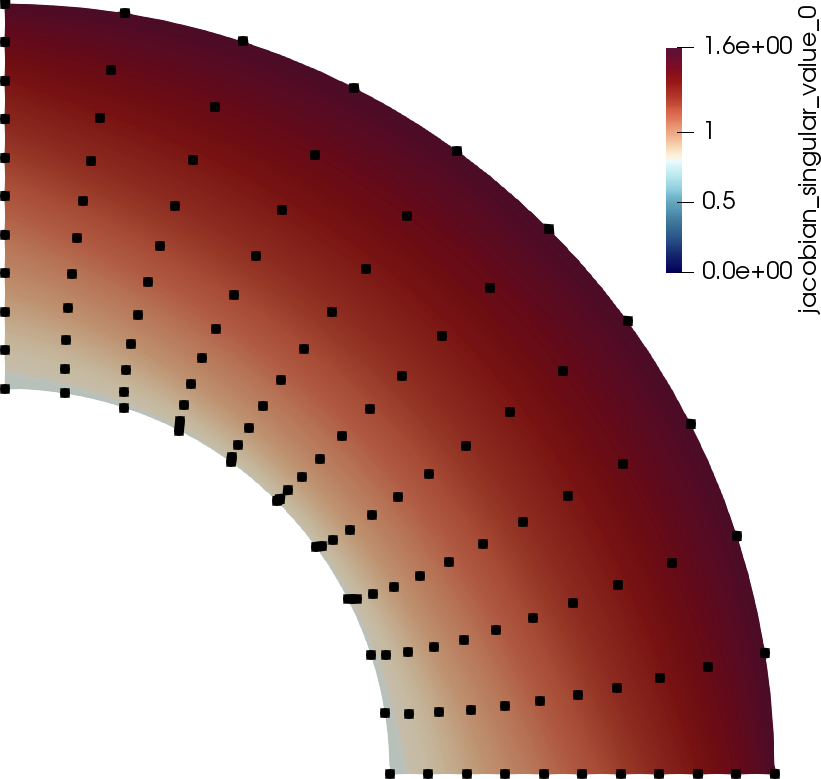}
    \includegraphics[width=0.22\textwidth]{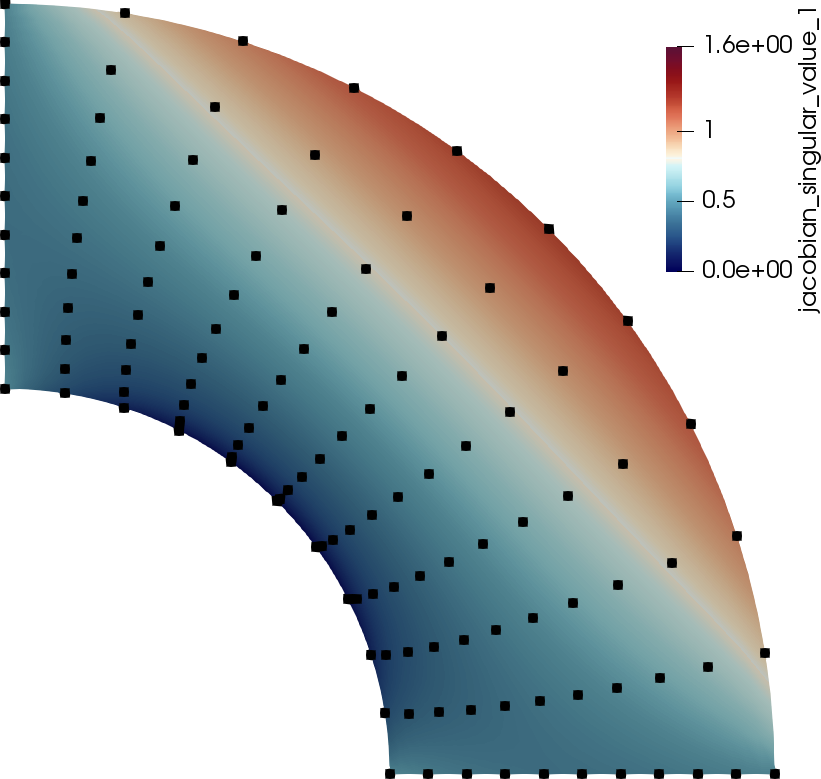}
    \hfill
    \includegraphics[width=0.22\textwidth]{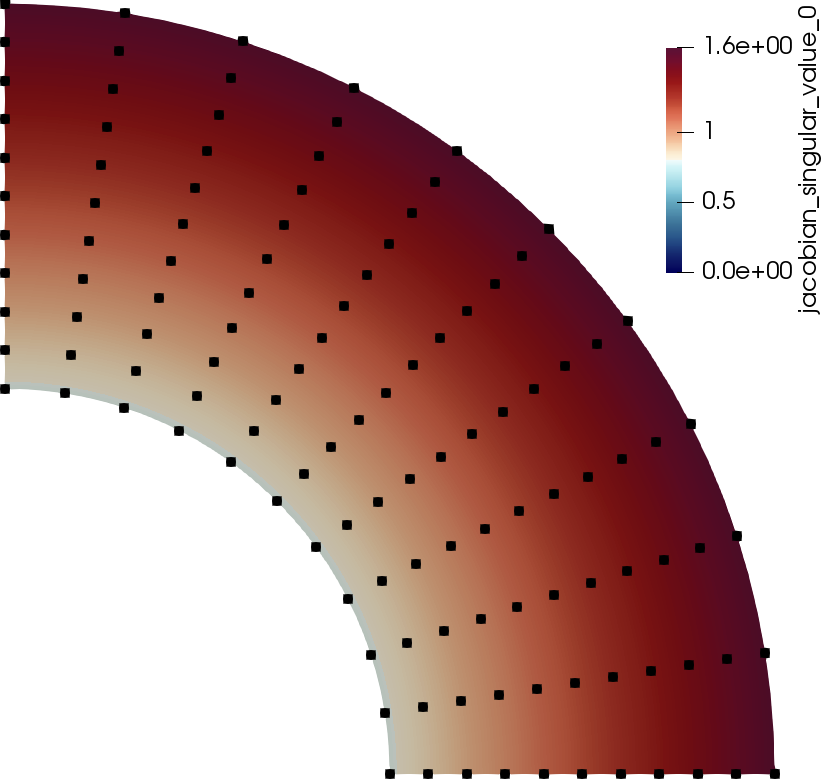}
    \includegraphics[width=0.22\textwidth]{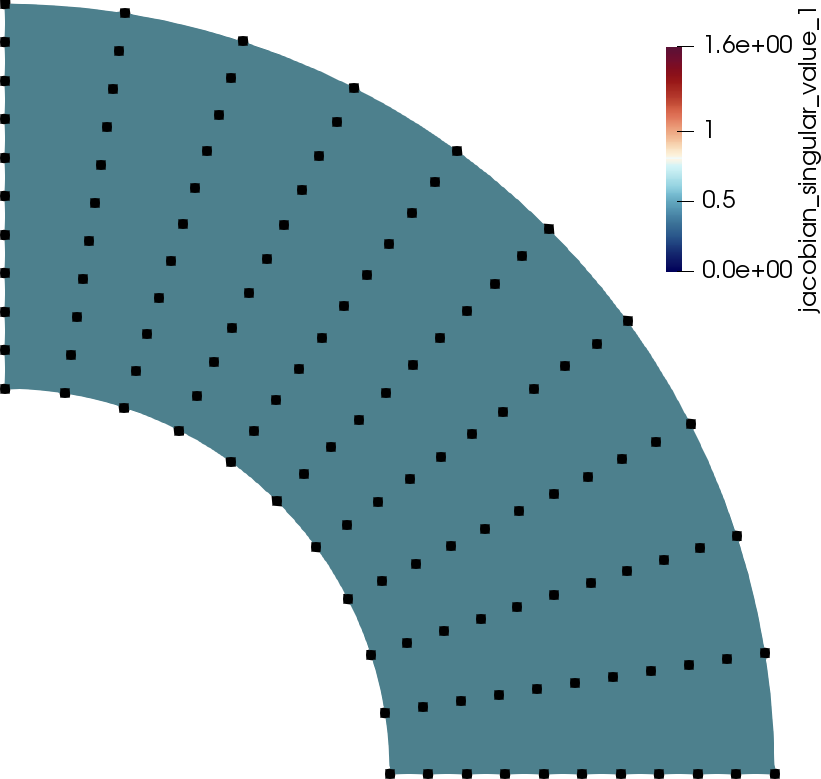}
    \caption{\it Illustration of the two singular values of $J_K=\hat\nabla \mathbf F_K(\mathbf {\hat x})$ for a quarter of a two-dimensional annulus; here, $\mathbf F_K$ is a polynomial mapping of a single element of degree 10. The ratio of the largest over the smallest singular
    value appears in the interpolation error estimate of the
    Bramble--Hilbert lemma, and consequently also in all
    error estimates for partial differential equations.
    Left two panels: Maximal and minimal singular value for point
    placement based on Laplace smoothing. Right two
    panels: Maximal and minimal singular value for point placement
    using a transfinite interpolation. (See the main text for interpretation.)
    The figures also include the positions of $11^2$ points
    equidistantly placed onto the reference cell and then mapped by
    $\mathbf F_K$ to the cell $K$ shown here, to illustrate the distortion.}
    \label{fig:transfinite}
\end{figure}


We visualize this approach using higher order
mappings. We recall that finite element error estimates on curved cells
depend on the product of a Sobolev
norm of higher order derivatives of $\mathbf F_K$
times a norm of the derivatives of $\mathbf F_K^{-1}$ (see, e.g., \cite{Ciarlet1972b} or
\cite[Sec. 3.3]{StrangFix88}), which we visualize by the ratio between the largest and
smallest singular value of $J_K$.
Figure~\ref{fig:transfinite} compares
the singular values for two variants of computing the interior points from
the surrounding 11 points per line, i.e., a mapping of polynomial degree
10, on a quarter of an annulus with inner and outer radii $0.5$ and $1$,
respectively. If the weights for the interior points are derived from
solving a Laplace equation in the reference coordinates, the representations
becomes distored, as is visible from the point distributions.
This leads to a ratio of up to 100 between the largest and smallest
singular value of the Jacobian $J_K=\hat\nabla \mathbf F_K(\mathbf {\hat x})$, and
theory suggests a break-down of convergence; in experiments,
$L_2$ errors of the solution to the Laplacian converge at best
at third order.
Conversely, using weights from transfinite interpolation results in a minimal
singular value of 0.5 throughout the whole domain in this example. The
resulting point distribution with transfinite interpolation in this particular
example is equivalent to an explicit polar description of the whole domain, but
applicable to generic situations with optimal convergence if
the coarse cells are valid. We note that these and similar concepts are
established in high-order
meshing, but with algorithms typically acting on the points of associated polynomial
descriptions, see e.g.~\cite{Solin04,Hindenlang2015,Moxey16,Mittal19,Dobrev19,Krais19,mengaldo2020industryrelevant} and references
therein, rather than the abstract definition used here.

We associate the transfinite interpolation with the
the initial (``coarse'') mesh of a finite element computation: each
coarse mesh cell is used to define the reference coordinate system
$\hat x$, and this is kept fixed even after many generations of
descendants. Interior edges between refined cells are then curved, ensuring high mesh quality,
assuming that the initial coarse cells reasonably approximate the
geometry. Some of the computations involved in this process can be
expensive, and we have therefore implemented caches that mitigate the cost for
the case of polynomial mappings.

\section{Application examples}
\label{sec:examples}

In the following, let us illustrate the ideas of the previous sections using concrete applications. In particular, we will show how the implementation of the two primitives affects the meshes one obtains for an industrial application (Section~\ref{sec:comparison-projection-cad}) and an example of how one can
choose metrics to generate graded meshes. In addition,
let us refer to Fig.~\ref{fig:hypershell} for an illustration of the transfinite
interpolation approach for extending surface descriptions to volume interiors.

\subsection{Surface meshes described by CAD geometries}
\label{sec:comparison-projection-cad}

As discussed in detail in previous sections, CAD surfaces consist of patches
that are the images of simpler domains in a two-dimensional $u$-$v$ space.
Let us first consider a case where the geometry is described by a single, albeit
rather complex, patch. The issue in even this
simplified case is that a patch can be parameterized in many
different ways, not all of which imply a more or less constant metric. 
Thus, it is unwise -- although very common -- to
generate a mesh in $u$-$v$ space to obtain the
corresponding three-dimensional surface grid, even though this of course has the advantage of generating nodes directly on the desired surface. Yet, the
resulting mesh will generally have cells of rather unequal sizes and may
show other severe deformations.

\begin{figure}
  \centering
  \begin{tabular}{c c c}
    \includegraphics[width=.23\textwidth]{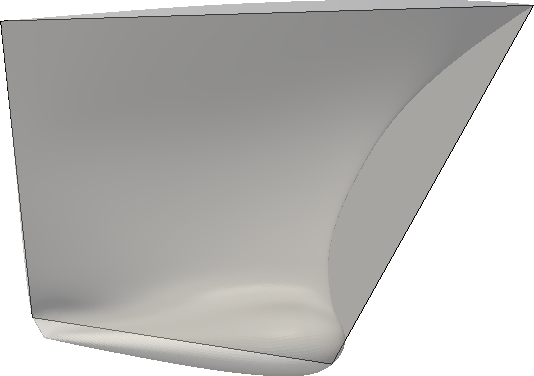} &
    \includegraphics[width=.23\textwidth]{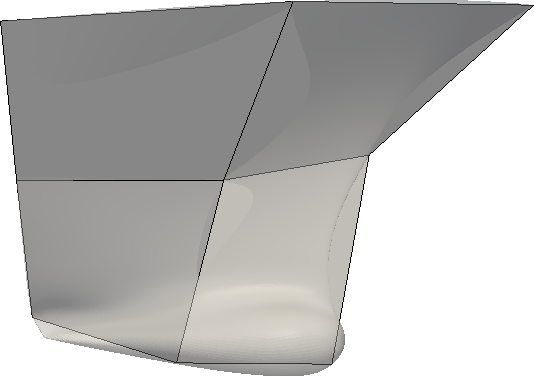} &
    \includegraphics[width=.23\textwidth]{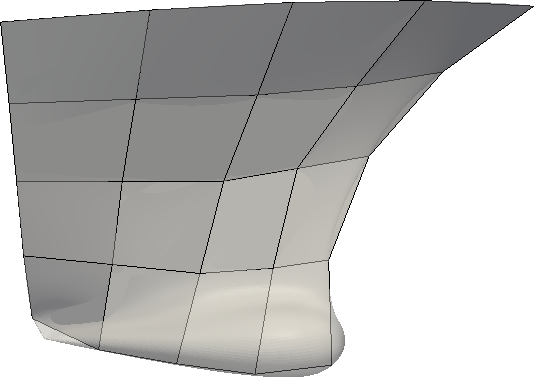} \\
    \includegraphics[width=.23\textwidth]{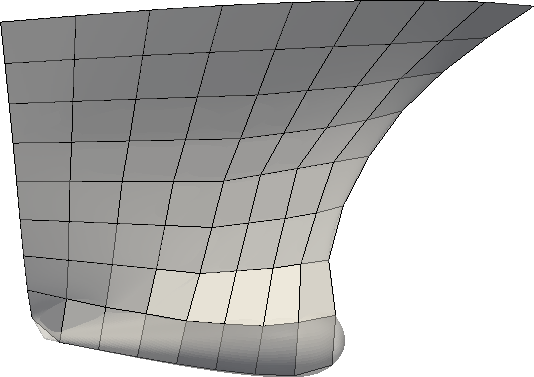} &
    \includegraphics[width=.23\textwidth]{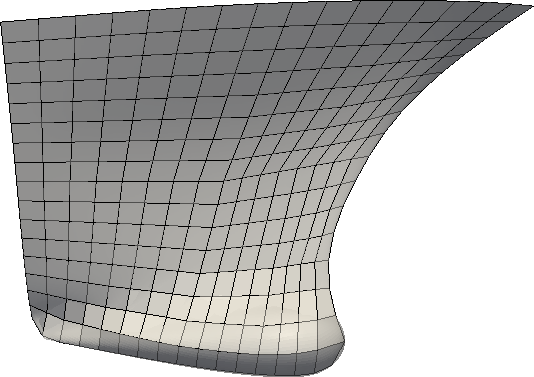} &
    \includegraphics[width=.23\textwidth]{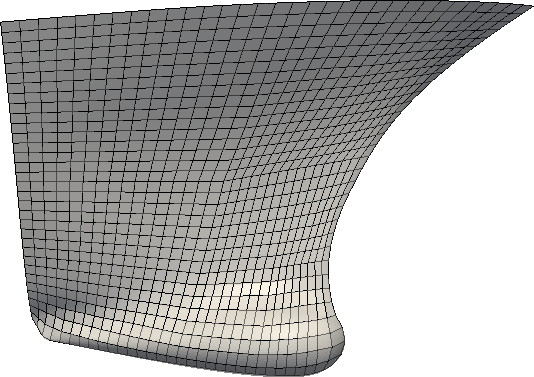} \\
    \includegraphics[width=.18\textwidth]{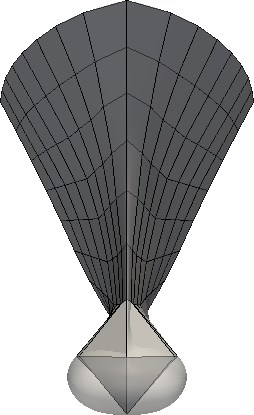} &
    \includegraphics[width=.18\textwidth]{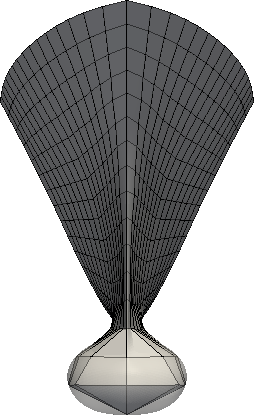} &
    \includegraphics[width=.18\textwidth]{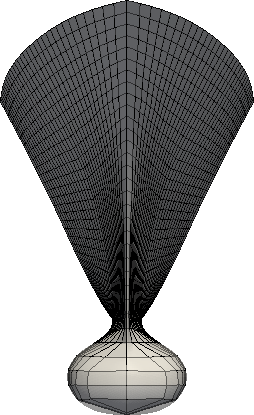}
  \end{tabular}
  \caption{\it Directional projection strategy with a horizontal
    direction of projection perpendicular to the axis of symmetry.
    The first two rows show
    side views of the coarse grid and grids obtained from five successive 
    refinements.  The last row shows a front view of
    the same grids shown in the second row. 
    This strategy produces uniformly distributed cells away from areas where the projection direction is close to the tangent to the shape
    (namely, at the bottom of the shape as well as the front of the
    bulb).}
  \label{fig:results-directional}
\end{figure}

\begin{figure}
  \centering
  \begin{tabular}{c c c}
    \includegraphics[width=.23\textwidth]{figures/step-54_common_0.png} &
    \includegraphics[width=.23\textwidth]{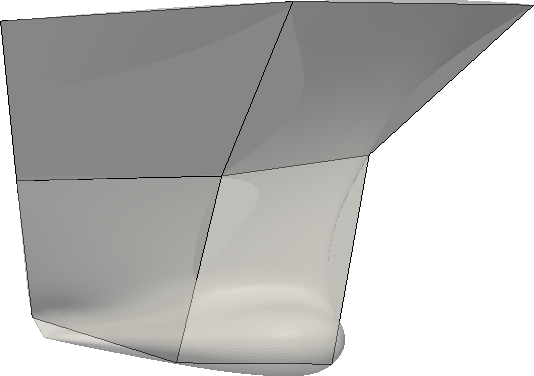} &
    \includegraphics[width=.23\textwidth]{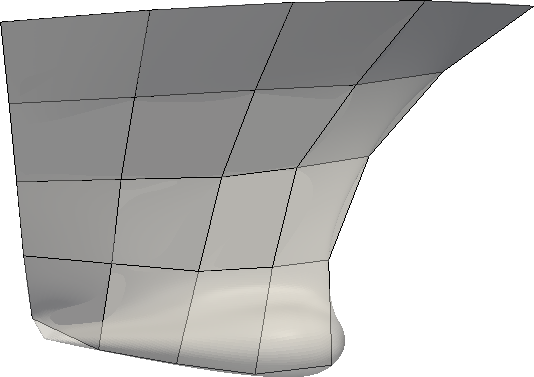} \\
    \includegraphics[width=.23\textwidth]{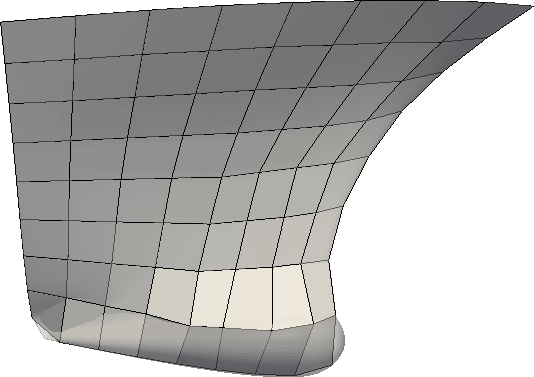} &
    \includegraphics[width=.23\textwidth]{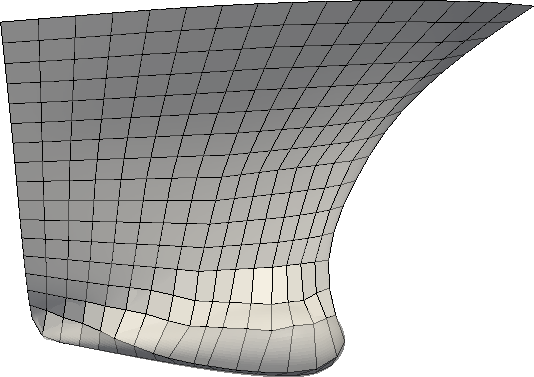} &
    \includegraphics[width=.23\textwidth]{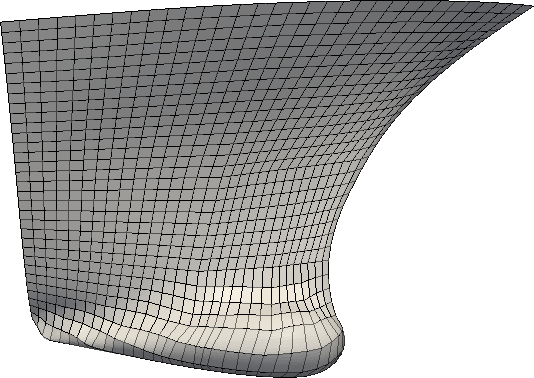} \\
    \includegraphics[width=.18\textwidth]{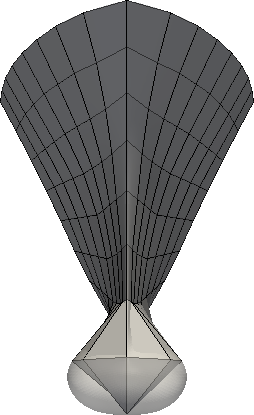}&
    \includegraphics[width=.18\textwidth]{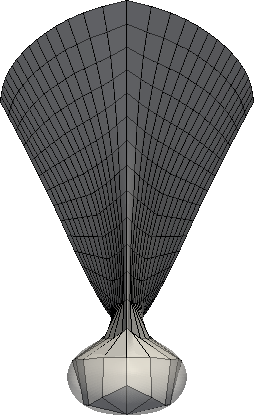} &
    \includegraphics[width=.18\textwidth]{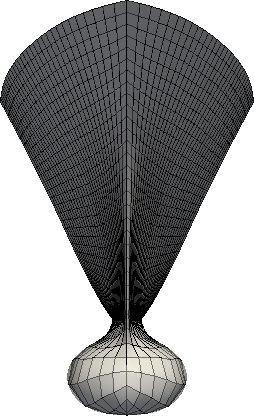}
  \end{tabular}
  \caption{\it Normal to mesh projection strategy. If the geometry
  does not intersect the direction normal to the existing points, then
  the closest point on the shape to the original point -- typically lying on the shape boundary -- is selected.
  Panels as in
  Fig.~\ref{fig:results-directional}.
    This strategy produces uniformly distributed cells in all cases.}
  \label{fig:results-normal-to-mesh}
\end{figure}

\begin{figure}
  \centering
  \begin{tabular}{c c c}
    \includegraphics[width=.23\textwidth]{figures/step-54_common_0.png} &
    \includegraphics[width=.23\textwidth]{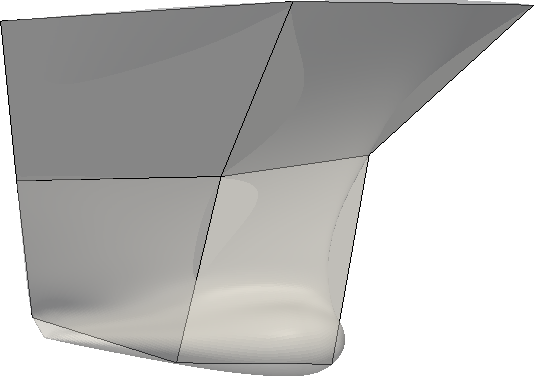} &
    \includegraphics[width=.23\textwidth]{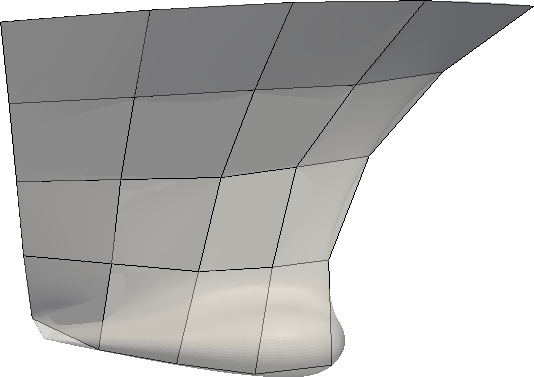} \\
    \includegraphics[width=.23\textwidth]{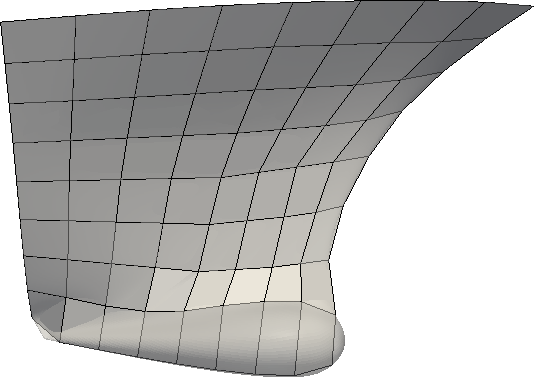} &
    \includegraphics[width=.23\textwidth]{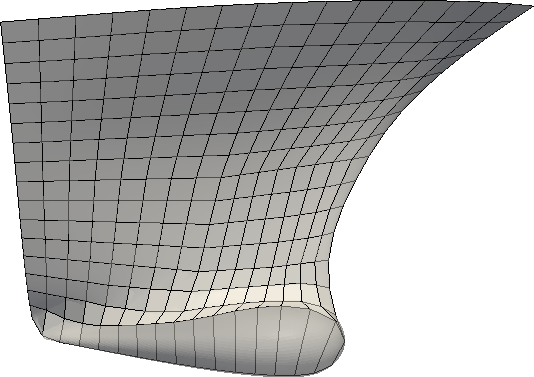} &
    \includegraphics[width=.23\textwidth]{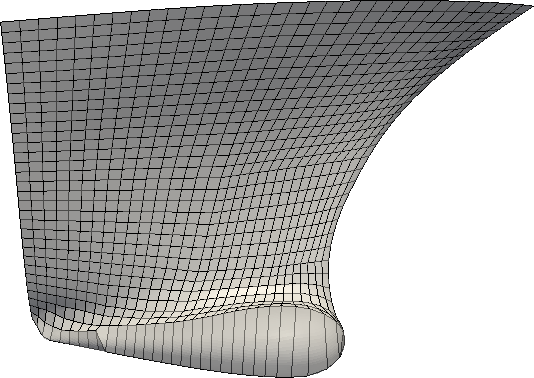} \\
    \includegraphics[width=.18\textwidth]{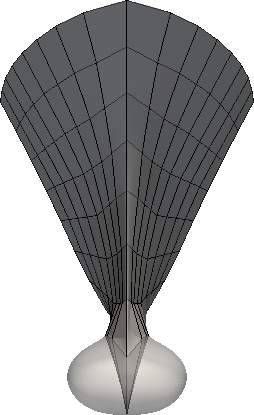} &
    \includegraphics[width=.18\textwidth]{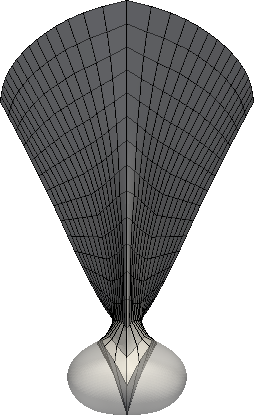} &
    \includegraphics[width=.18\textwidth]{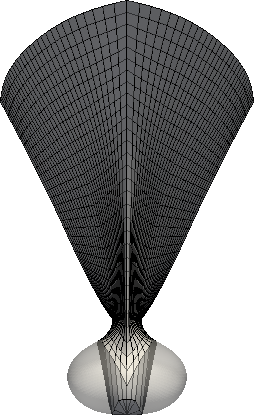}
  \end{tabular}
  \caption{\it Normal to surface projection strategy. In cases where more
  than one surface normal projection is available, the closest of them is
  selected. If the shape is composed by several sub-shapes, the
  projection is carried out onto every sub-shape and the closest
  projection point is selected.
  Panels as in
  Fig.~\ref{fig:results-directional}.
    This strategy is unable to produce well-shaped cells in areas of large curvature.}
  \label{fig:results-normal}
\end{figure}

To illustrate how our approach can be used to generate better meshes, we will use an industrial application that
involves meshing a single parametric patch
describing the bow portion of one side of
the DTMB 5410 ship hull, containing also a sonar dome. The presence of several
convex and concave high curvature regions makes such a geometry a
particularly meaningful example.

Figures~\ref{fig:results-directional}--\ref{fig:results-normal} show
results for this model geometry with the three projection strategies.
The directional projection strategy with a horizontal
direction of projection (Fig.~\ref{fig:results-directional})
generally produces high quality meshes except
in those places where the geometry is tangent to the
projection direction -- i.e., in particular at the front of the bulb
as well as the bottom. In contrast, using the direction normal to
the existing points (Fig.~\ref{fig:results-normal-to-mesh})
generates high quality meshes everywhere. Finally, the option to
use a surface normal instead of a mesh normal vector
(Fig.~\ref{fig:results-normal}) is not only expensive to compute,
but here yields meshes that are grossly distorted
wherever the geometry has large curvature (i.e., around the bulb); this
might also have been expected from the left panel of
Fig.~\ref{fig:surface-projection-comparison} that shows a similar effect.

\subsection{Refinement strategy based on local maximum curvature}
\label{sec:local-curvature-example}

The case of multi-patch geometries is more complicated as
small gaps or superimpositions are typically present between neighboring
patches. Meshes directly generated from this parameterization will therefore
often not be ``water-tight''. On the other hand, one can generate a coarse mesh
by hand or by software that simply starts with a few points on the surface
that are then connected to cells without taking into account the subdivision
into patches; such a mesh can then be refined hierarchically to obtain
a mesh of sufficient density.

\begin{figure}
  \phantom{.}
  \hfill
  \includegraphics[height=3cm]{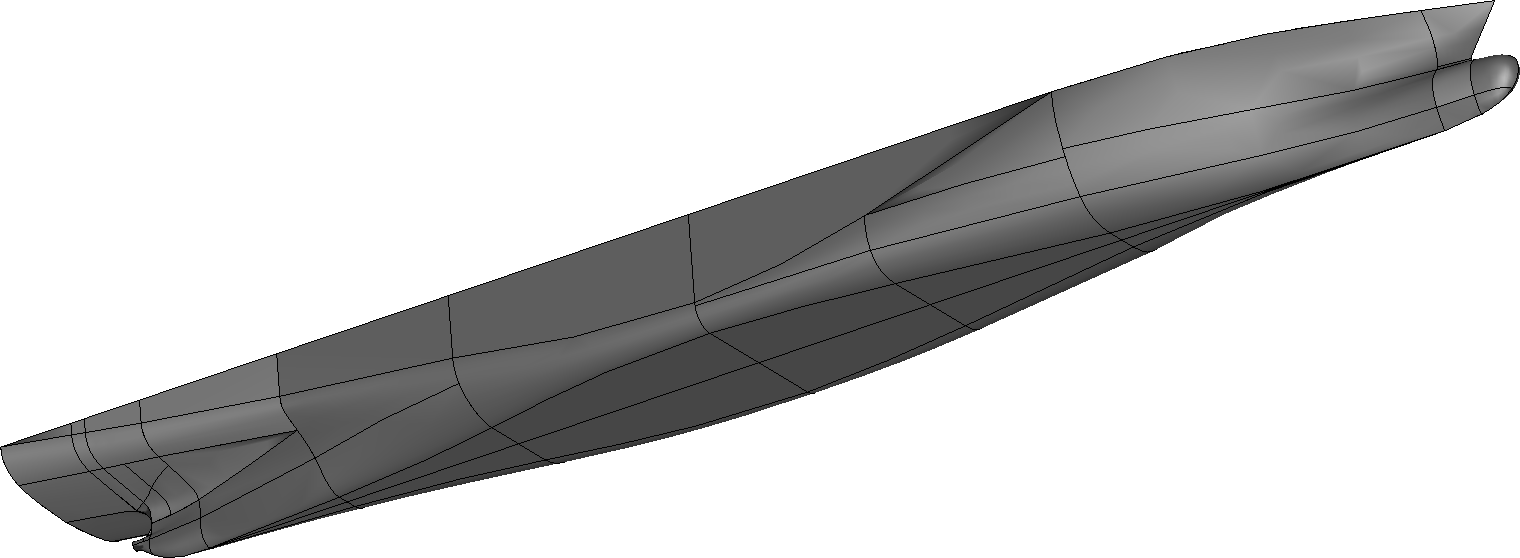}
  \hfill
  \includegraphics[height=3cm]{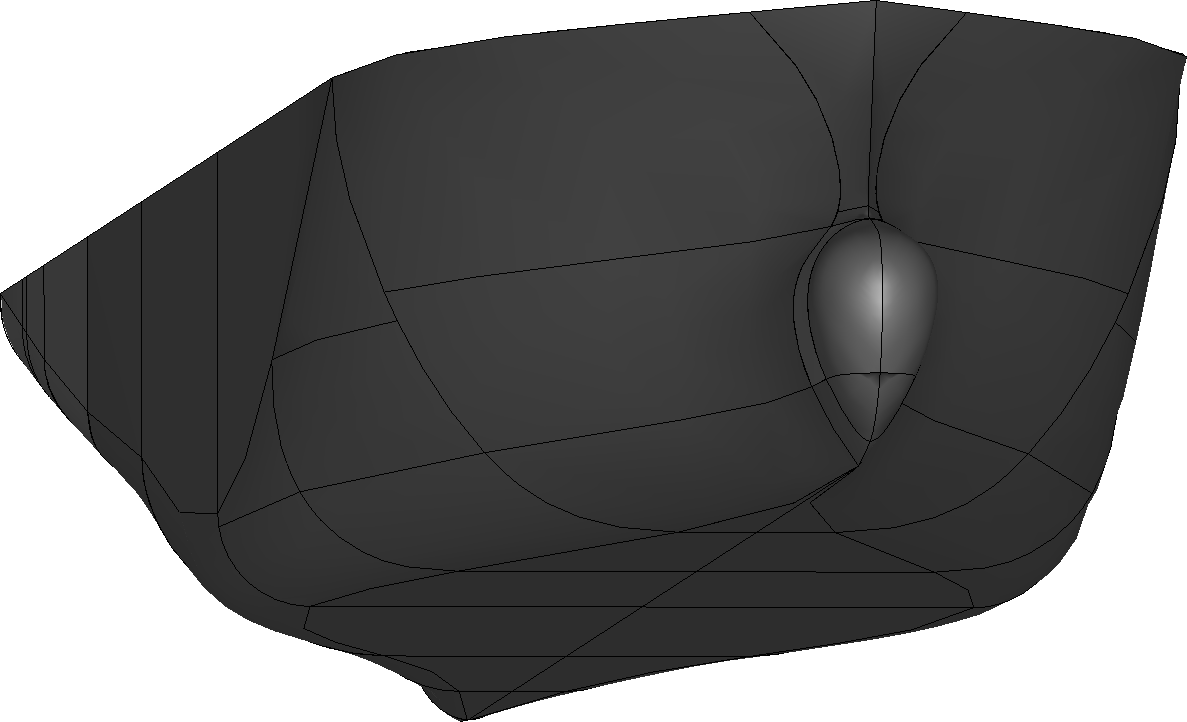}
  \hfill
  \phantom{.}
  \caption{\it Overall (left) and bow (right) view of the CAD model of a
           ship hull. The model is composed of approximately 120 parametric patches, delimited by black lines.}
  \label{fig:kcs_cad}
\end{figure}

In the following, let us show how we can use the results of the previous section towards building meshes
for an entire ship hull; we will also show how additional information can be extracted from CAD tools to drive the refinement strategy on CAD based geometries.
To this end, Figure \ref{fig:kcs_cad} depicts a CAD model of a Kriso KCS ship hull
-- a common benchmark for CFD applications of naval architecture \cite{Kim_et_al_Kriso_2001}. In this
production-like CAD model the patches are not connected in a water-tight fashion and the
surface parametrization is not continuous at patch junctions.
These defects prevent most mesh generators from obtaining a closed grid. 

\begin{figure}
\centering
\includegraphics[width=\textwidth]{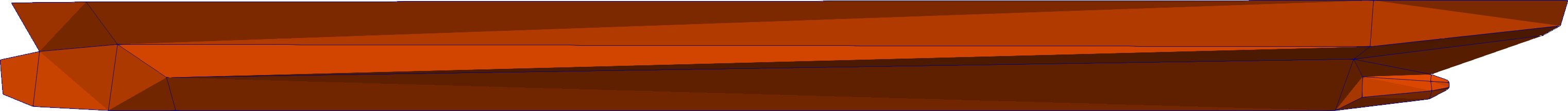}
\\[.5cm]
\includegraphics[width=\textwidth]{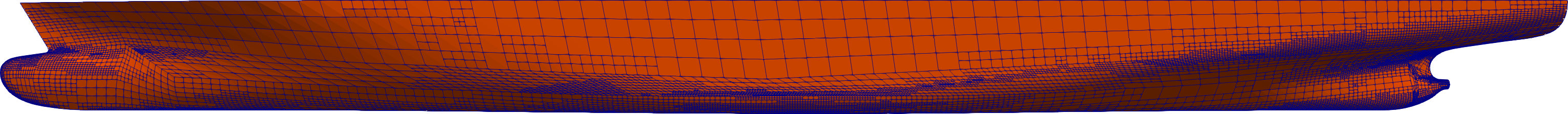}
\caption{\it Initial (top) and final (bottom) surface grids on the
         Kriso KCS hull, with 40 and $\tilde{}$
         11,500 quadrilateral cells, respectively.}
  \label{fig:kcs_grid}
\end{figure}

We start from a minimal initial surface grid composed of about 40 cells 
(top panel of Fig.~\ref{fig:kcs_grid}) and refine it a number of times
using the strategy where we project in a direction normal to the existing points.
Refinement of the initial grid starts with an anisotropic refinement step in which
cells with an aspect ratio larger then a threshold $\lambda_\text{max}$ are cut
along their most elongated direction. 
We then refine the resulting quadrilateral mesh adaptively, using an estimate of the
local curvature of the CAD
surface as a criterion. This estimate is obtained exploiting a technique essentially identical to the one used in the ``Kelly'' error estimator~\cite{KGZB83,GKZB83}. Elements in areas with higher curvature will have larger jumps of the (cell) normal
vectors across cell boundaries. 
The bottom panel of Fig.~\ref{fig:kcs_grid} shows the final grid generated. 

\begin{figure}
\centering
\includegraphics[width=0.49\textwidth]{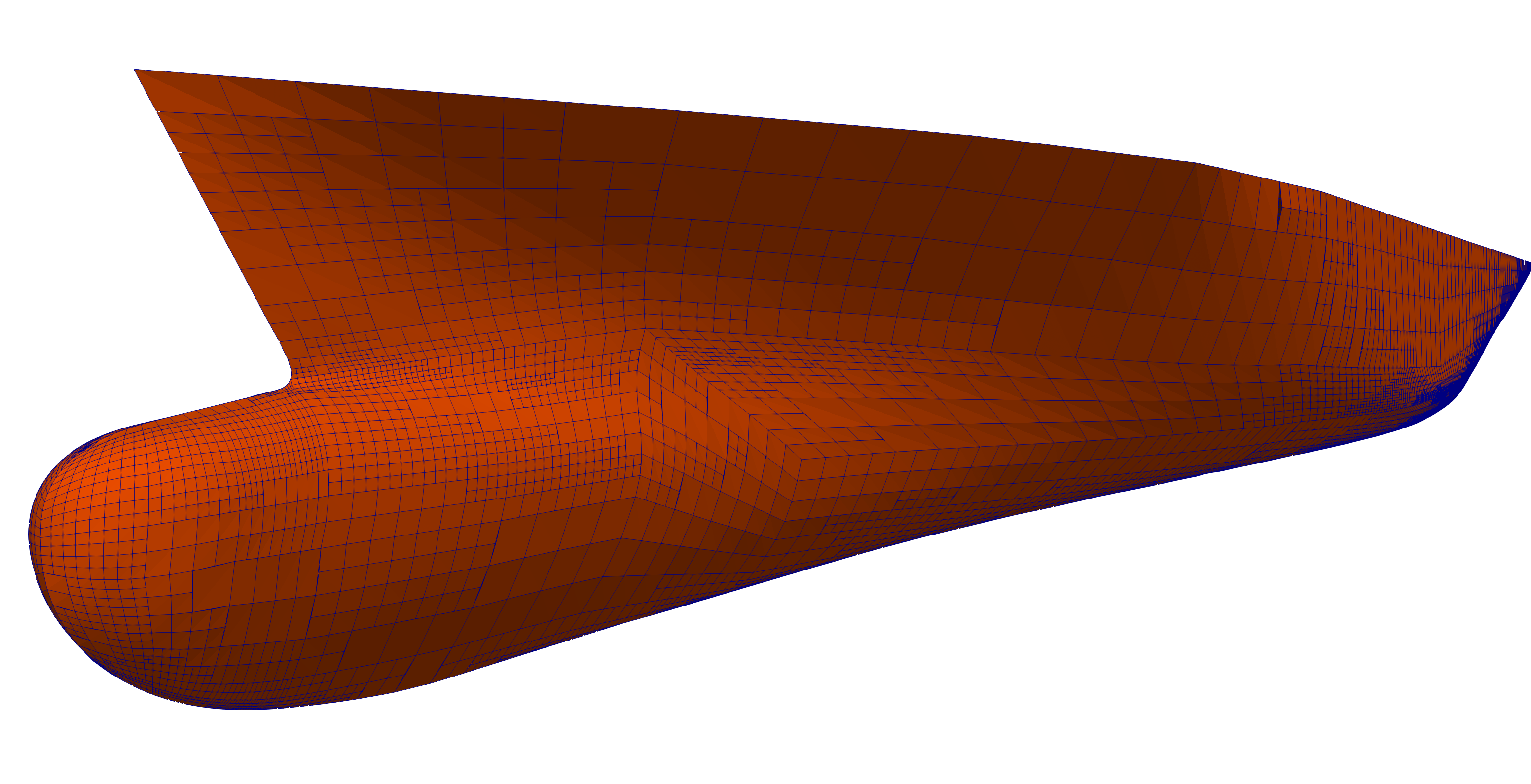}
\hfill
\includegraphics[width=0.49\textwidth]{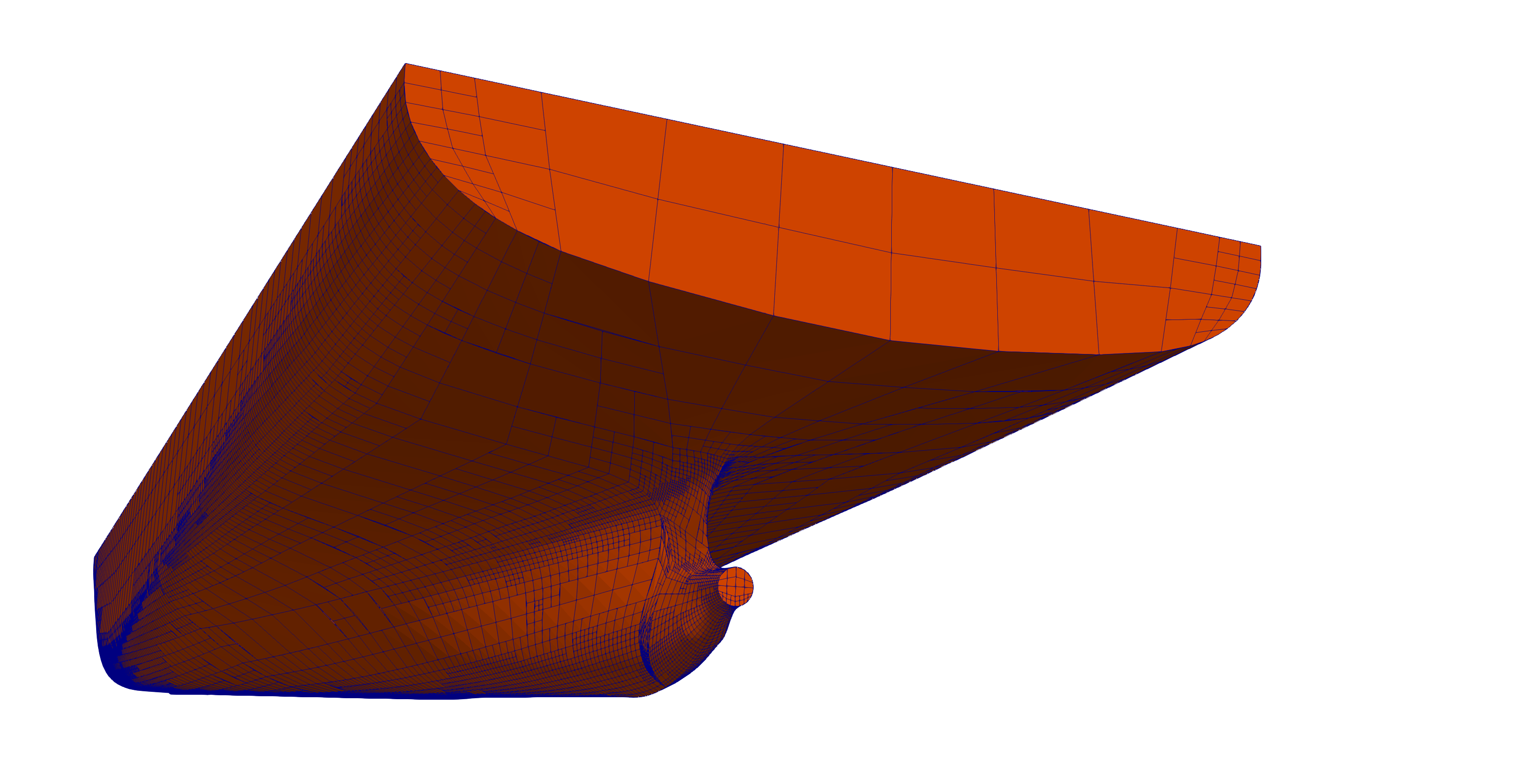}
\caption{\it Bow (left) and stern (right) details of the adaptively refined
         surface grid on the Kriso KCS hull from Figures \ref{fig:kcs_cad} and \ref{fig:kcs_grid}. 
         The adaptive refinement strategy results in finer cells in high
         curvature regions, ensuring uniform approximation of the geometry.
         In addition, the grid is independent
         of the non-sharp edges separating the 120 parametric patches composing
         the underlying CAD model. In the final mesh,
         hanging nodes are placed on the underlying geometry, leading to a non-watertight mesh. However, these 
         artifacts are easily removed by enforcing continuity 
         of the geometry.}
  \label{fig:kcs_grid_detail}
\end{figure}

Despite the fact that the original CAD surface is composed of
several unconnected parametric patches, the projection procedure is able to find,
for each of the grid nodes, the best projection among those obtained onto individual patches.
As a result, the final mesh is water-tight, and
independent of the CAD surface parametrization and patch distribution. In
addition, the adopted refinement strategy distributes a larger number of new nodes
in high curvature regions, ensuring a uniform quality of approximation of the
geometry. We show this in more detail in
Fig.~\ref{fig:kcs_grid_detail}, illustrating the bow and stern
portions of the final grid. Finally, it is worth pointing out that for the generation of the
grids portrayed in Fig.~\ref{fig:kcs_grid} and Fig.~\ref{fig:kcs_grid_detail}, no smoothing
stage was carried out in between refinements to enhance mesh quality. Yet, the projection in a direction normal to the existing points allows
for retaining the quality of the original coarse grid across more than 10 levels of refinement without any additional adjustment.

\subsection{(Ab)Using primitives based on analytic charts for graded meshes}
\label{sec:abusing}

We end this section by outlining how else the primitives can be used
to satisfy practical needs. In particular, we can achieve graded mesh refinement strategies based on the explicit definition of a custom metric that describes the manifold with associated push-forward and pull-back operations as discussed in Section~\ref{sec:chart-case}.

\begin{figure}
\centering
  \phantom{.}
  \hfill
\includegraphics[width=.22\textwidth]{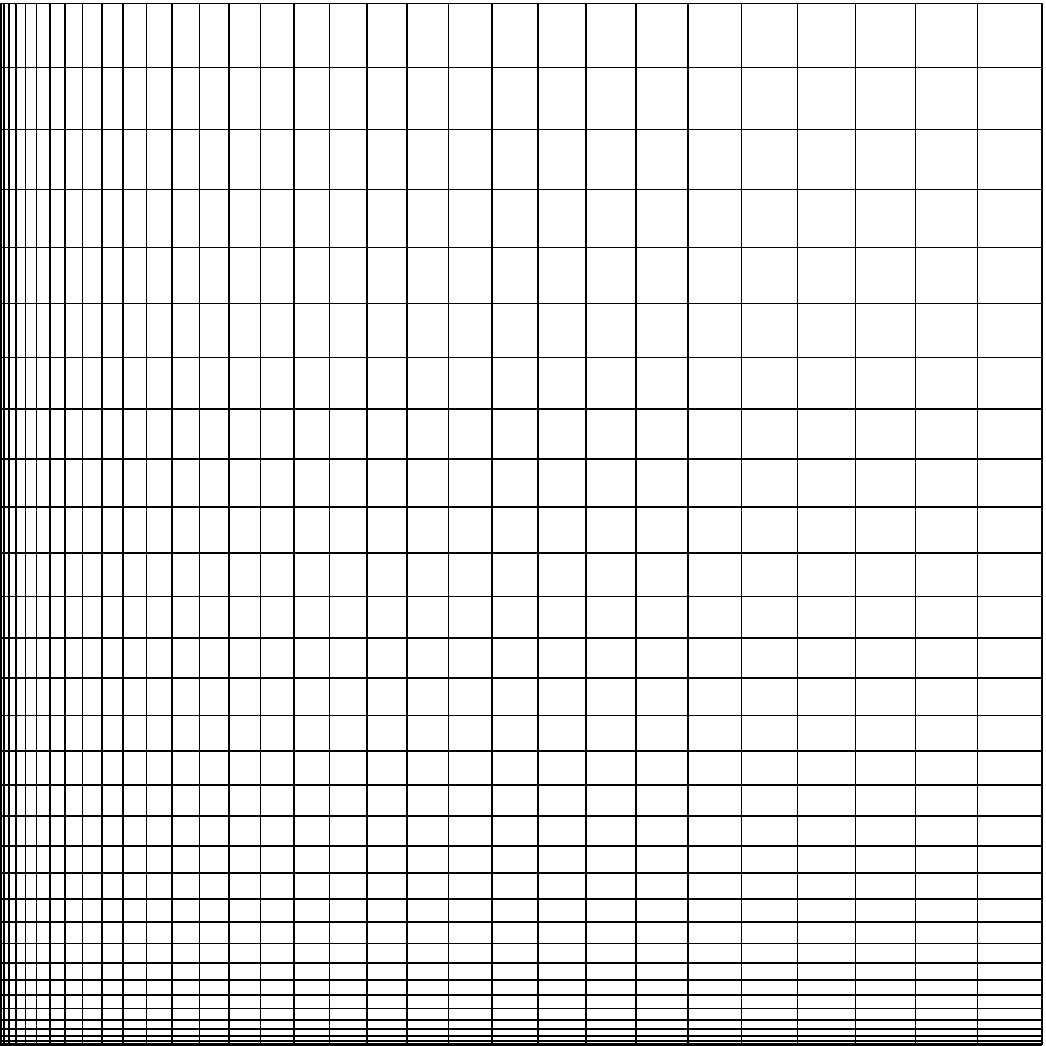}
	\hfill
    \includegraphics[width=.22\textwidth]{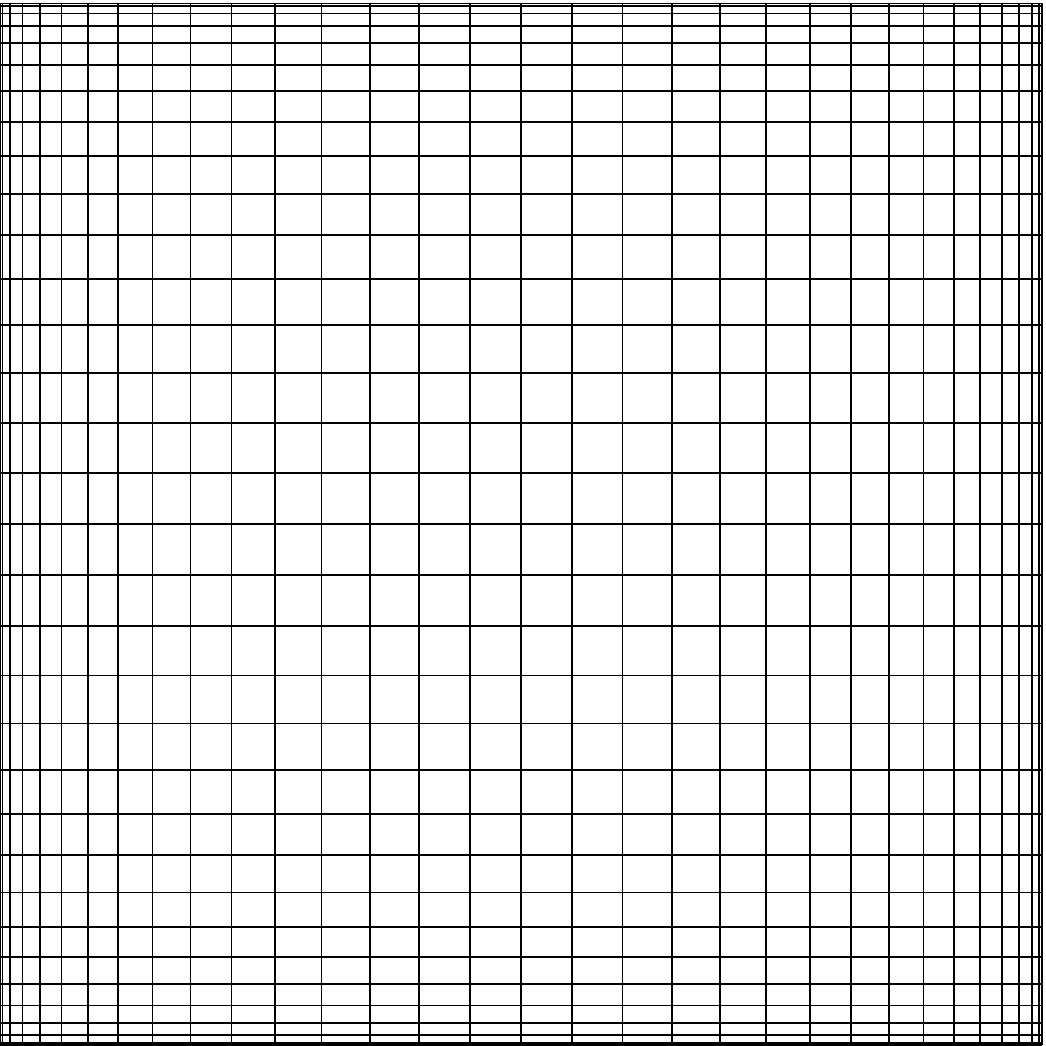}
    \hfill
    \phantom{.}
\caption{\it Example of graded meshes, obtained using simple diffeomorphisms of the unit square.}
\label{fig:graded-mesh}
\end{figure}

Consider, as a simple example, the discretization of the unit square $[0,1]^2$. We can provide a non-flat local metric, induced by the following (analytic, invertible) mapping:
\begin{equation*}
 \begin{split}
		\phi: [0,1]^2 & \to [0,1]^2 \\
        (x, y) &\to (x^2, y^2).
 \end{split}
\end{equation*}
If we consider this mapping for the construction of the \textsc{new point} primitive as described in Section~\ref{sec:chart-case}, we obtain, upon refinement, a graded mesh like in Figure~\ref{fig:graded-mesh}, left.
But we are not restricted to refining towards left and bottom edges in the figure. Rather, by using more complicated mappings, one can
 concentrate elements in regions of interest by choosing $\phi$. For 
 example, we can refine towards all four sides of the
 square (see  Figure~\ref{fig:graded-mesh}, right) by choosing
\begin{equation*}
		\phi(x, y) = \left(\frac12\sin\left(\pi\left(x-\frac12\right)\right)+\frac12, \frac12\sin\left(\pi\left(y-\frac12\right)\right)+\frac12\right).
\end{equation*}

\section{Conclusions}
\label{sec:conclusions}

The traditional workflow in finite element, finite volume, or finite difference simulations -- split into preprocessing, simulation, and
postprocessing stages -- is poorly suited to modern numerical methods
because information about the underlying geometry is typically not
propagated beyond the mesh generation stage. To address this deficiency,
we have herein undertaken a comprehensive review
of where geometry information is used in simulation software, and
what it would take to provide this from typical Constructive Solid
Geometry (CSG) or Computer Aided Design (CAD) tools. We have found
that, despite the very large number of places in which geometry
information is used, all uses can be reduced to just two ``primitive''
operations: The generation of a new point that interpolates existing
points, and the computation of a tangent vector to a line connecting
two points. We have then described in detail how these two operations
can be implemented for common CSG and CAD cases, and illustrated how
these implementations can help create meshes for complex geometries
that are both of high quality and respect the underlying geometry.

The methods discussed herein are all available as part of the
open source finite element library deal.II \cite{dealiicanonical,dealII91}
and form the basis of propagating geometry information to all
parts of the library and application codes built on top of it.

\appendix

\section{An operation that can not be expressed via the two primitives}
\label{sec:exceptions}

In our survey of uses of geometric information for Section~\ref{sec:uses}, we have come across 
just one operation that can not be expressed in terms of
the two primitives introduced in Section~\ref{sec:primitives}:
the computation of second or higher derivatives of finite
element shape functions when using the \textit{underlying} geometry of the
domain, rather than a polynomial mapping.

For most ``common'' finite elements, shape functions are defined on a
reference cell and then mapped to each of the cells of the finite
element mesh. Let $\hat K$ be the reference cell (e.g., the reference
triangle or reference square) and $\hat\varphi_i(\hat{\mathbf x}),
i=1,\ldots,N$ be the shape functions of a scalar finite element. For a
given cell $K$ of the mesh, let $\mathbf F_K:\hat K \mapsto K$ be the mapping,
assuming that $K$ is part of the underlying geometry,
rather than a polynomial approximation of it. In that case, $\mathbf F_K$ is,
in general, not polynomial. The shape functions we will then work with
are $\varphi_i(\mathbf x) = \hat\varphi_i(\mathbf F_K^{-1}(\mathbf
x))$.

By the chain rule, the derivatives, in real space, of shape
functions at an evaluation point $\hat{\mathbf x}_q$ are then given by
\begin{align*}
  \nabla\varphi_i(\mathbf x_q) 
  &= J_K^{-1} \hat\nabla\hat\varphi_i(\hat{\mathbf x}_q),
\end{align*}
where $J_K=\frac{\partial \mathbf x}{\partial \hat{\mathbf
    x}}=\frac{\partial \mathbf F_K(\hat{\mathbf x})}{\partial
  \hat{\mathbf x}}$ is the Jacobian of the transformation. This
matrix, and its inverse, are easily computed for both polynomial
mappings (see Section~\ref{sec:mappings}) and for exact mappings (see
Section~\ref{sec:jacobian}). However, for some applications,
it is also necessary to compute second or even higher derivatives of
shape functions. An example is the evaluation of the residual of a
finite element solution for second or higher order operators.
Others include the evaluation of the bilinear form for the biharmonic
equation, and the evaluation of the Hessian of (a component of) the
solution for purposes of defining anisotropic mesh refinement indicators. In these
cases, we can compute (with appropriate contraction over the various
indices of the objects in the formula)
\begin{align*}
  \nabla^2 \varphi_i(\mathbf x_q) 
  &= J_K^{-1} \hat\nabla\left[
    J_K^{-1} \hat\nabla
    \hat\varphi_i(\hat{\mathbf x}_q)
    \right]
  \\
  &= J_K^{-1} J_K^{-1} \hat\nabla^2 \hat\varphi_i(\hat{\mathbf x}_q)
     + J_K^{-1} \left[ \hat\nabla [J_K^{-1}] \right]
    \hat\varphi_i(\hat{\mathbf x}_q)
  \\
  &= J_K^{-1} J_K^{-1} \hat\nabla^2 \hat\varphi_i(\hat{\mathbf x}_q)
     - J_K^{-1} \left[ J_K^{-1} (\hat\nabla J_K) J_K^{-1}\right]
    \hat\varphi_i(\hat{\mathbf x}_q),
\end{align*}
where we have made use of the equality $0=\hat\nabla
(J_KJ_K^{-1}) = (\hat\nabla J_K)J_K^{-1} + J_K(\hat\nabla [J_K^{-1}])$.

This representation of second derivatives (and formulas
for even higher derivatives) requires computing derivatives of
$J_K$. This is not a complication for the usual, polynomial
mappings such as those discussed in Section~\ref{sec:mappings}, since
all we need to do is determine the appropriate interpolation points,
construct the polynomial approximation, and then take derivatives of
this polynomial. On the other hand, for ``exact'' mappings, no
explicit representation of the mapping is available.
We can only use the \textsc{new point} primitive to evaluate
$\mathbf F_K$, and use the \textsc{tangent vector} primitive to
compute the derivative $J_K=\hat\nabla \mathbf F_K$, but we have no way
of evaluating the derivative of $J_K$ (i.e., the second derivatives of
$\mathbf F_K$). 

\begin{remark}
  \label{remark:third-primitive}
  One could of course define a third primitive to also aid the
  computation of this information. (See also Remark~\ref{remark:minimality}.)
  On the other hand, our goal in this
  paper is to point out a \textit{minimal} interface that geometry packages have to
  satisfy to support {\em the most common} finite element operations. Since
  all finite element packages we are aware of support polynomial geometry
  approximations, we consider the two operations defined in
  Section~\ref{sec:primitives} as sufficient and defer the definition
  of additional primitives necessary for the operation discussed in
  this section to separate work.

  Alternatively, one can approximate $\hat\nabla J_K$ by finite
  differencing the matrix $J_K$, which we already know how to
  compute.
\end{remark}

\section{Periodic domains}
\label{sec:periodicity}

In practice, one frequently encounters domains that are
periodic in at least one direction, or for which  at least the
\textit{push-forward} function is periodic. An example is a domain
that consists of the surface of a cylinder. Another would be a
situation where one only considers a segment of a cylinder surface --
in that case, the domain itself is not periodic, but both the
\textit{image} and pre-image of the push-forward function are. Such
cases present interesting questions for the implementation of the
\textsc{new point} primitive. If we can answer these appropriately, we
will know how to deal with the \textsc{tangent vector} primitive by using
relationship~\eqref{eq:tangent}.

\begin{figure}
    \centering
    \includegraphics[width=0.6 \textwidth]{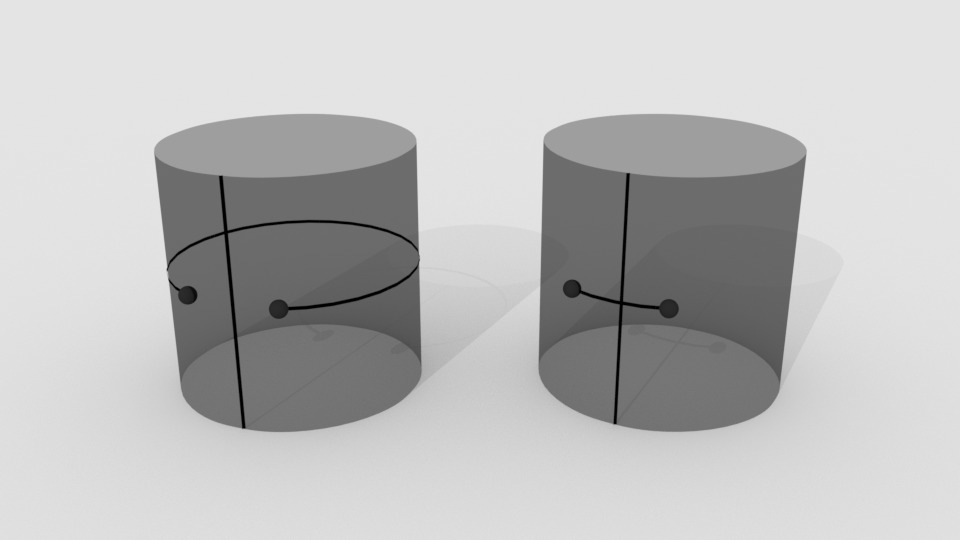}
    \caption{\it Effect of periodicity when computing the \textsc{new point} primitive with two points close to the gluing edge of a cylinder. If periodicity is not taken into account, the path used to compute the averages does not cross the gluing edge (left), and new points would end up on the line that turns around the cylinder. When taking periodicity into account, the shortest path is used when computing averages (right), and new points end up on the line that crosses the gluing edge.}
    \label{fig:periodic-domain}
\end{figure}

In our implementation in \dealii{}, primitives based on pull-back and push-foward functions are provided by the \texttt{ChartManifold} class
that allows one to specify along which direction in the reference Euclidean patch $U_\alpha^\ast$ periodicity should be considered, and what the periodicity period $L$ along that direction is. Specifically, periodicity affects the way the \textsc{new point} oracle computes the middle point between two neighboring points: If two points are more than half a period distant in the Euclidean reference patch $U_\alpha^\ast$, then their distance should be computed by crossing the periodicity boundary, i.e., the weighted average of the two points should be computed by adding a half period to the sum of the two. For example, if along a direction  we have a periodic manifold, with period of $2\pi$, then the average of $(2\pi-\epsilon)$ and $(\epsilon)$ along the periodic domain (computed with equal weights) should not return $\pi$, but $2\pi$ (or, equivalently, zero). On the manifold, these two points are truly at a distance $2\epsilon$ and not $(2\pi-\epsilon)$ (see Figure~\ref{fig:periodic-domain} for an example of this process).

The periodicity treatment becomes ill posed for cases in which the
distance of two points is \textit{exactly} a half period. Then, either
direction would be possible, and no unique solution exists. The straight forward way to solve this issue is to ensure that there are no points at exactly this distance, by providing a slightly more refined initial grid.

\section*{Acknowledgements}
L.~Heltai and A.~Mola were partially supported by the PRIN grant No.~201752HKH8, ``Numerical Analysis for Full and Reduced Order Methods for the efficient and accurate solution of complex systems governed by Partial Differential Equations (NA-FROM-PDEs)''. A.~Mola was also partially supported by the project UBE2 -
``Underwater blue efficiency 2'' funded by Regione FVG, POR-FESR 2014-2020, Piano Operativo Regionale Fondo Europeo per lo Sviluppo Regionale.
  W.~Bangerth was partially
supported by the National Science Foundation under award OAC-1835673
as part of the Cyberinfrastructure for Sustained Scientific Innovation (CSSI)
program; by award DMS-1821210; 
by award EAR-1925595;
and by the Computational Infrastructure
in Geodynamics initiative (CIG), through the National Science
Foundation under Awards No.~EAR-0949446 and EAR-1550901 and The
University of California -- Davis.
 M.~Kronbichler was partially supported by the German Research Foundation (DFG) under
    the project ``High-order discontinuous Galerkin for the exa-scale''
    (ExaDG) within the priority program ``Software for Exascale Computing''
    (SPPEXA), grant agreement no.~KR4661/2-1.

\bibliographystyle{abbrv}
\bibliography{paper}

\end{document}